\documentclass[11pt]{amsart}

\usepackage[all]{xy}
\usepackage{amscd}
\usepackage{mathrsfs}
\usepackage{amssymb}
\usepackage{amsthm}
\usepackage{amsmath}

\newtheorem{thm}{Theorem}[section]
\newtheorem{prop}[thm]{Proposition}
\newtheorem{lem}[thm]{Lemma}
\newtheorem{cor}[thm]{Corollary}
\newtheorem{rem}[thm]{Remark}

\newtheorem{ex}[thm]{Example}
\newtheorem{exercise}[thm]{Exercise}
\newtheorem*{ithm}{Theorem A (Local Class Field Theory)}

\theoremstyle{definition}
\newtheorem{defn}[thm]{Definition}

\newcommand{\pfbegin}{\begin{proof}}
\newcommand{\pfend}{\end{proof}}

\newcommand{\bde}{\begin{defn}}
\newcommand{\ede}{\end{defn}}
\newcommand{\bpr}{\begin{prop}}
\newcommand{\epr}{\end{prop}}
\newcommand{\bth}{\begin{thm}}
\newcommand{\ethm}{\end{thm}}
\newcommand{\bexa}{\begin{ex}}
\newcommand{\eexa}{\end{ex}}
\newcommand{\bex}{\begin{exercise}}
\newcommand{\eex}{\end{exercise}}
\newcommand{\bcor}{\begin{cor}}
\newcommand{\ecor}{\end{cor}}
\newcommand{\blem}{\begin{lem}}
\newcommand{\elem}{\end{lem}}
\newcommand{\brem}{\begin{rem}}
\newcommand{\erem}{\end{rem}}
\newcommand{\bprf}{\begin{proof}}
\newcommand{\eprf}{\end{proof}}

\begin{document}

\newcommand{\ol}{\overline}
\newcommand{\ra}{\rightarrow}
\newcommand{\Lra}{\Longrightarrow}
\newcommand{\lra}{\longrightarrow}
\newcommand{\N}{\mathbb{N}}
\newcommand{\Z}{\mathbb{Z}}
\newcommand{\Q}{\mathbb{Q}}
\newcommand{\R}{\mathbb{R}}
\newcommand{\F}{\mathbb{F}}
\newcommand{\G}{\mathbb{G}}
\newcommand{\CO}{\mathcal O}
\newcommand{\mO}{\mathcal O}
\newcommand{\fkp}{\mathfrak p}
\newcommand{\fkm}{\mathfrak m}
\newcommand{\Gal}{\operatorname{Gal}}
\newcommand{\gal}{\operatorname{gal}}
\newcommand{\ke}{\operatorname{Ker}}
\newcommand{\im}{\operatorname{Im}}
\newcommand{\Frob}{\operatorname{Frob}}
\newcommand{\Fr}{\operatorname{Fr}}
\newcommand{\Image}{\operatorname{Im}}
\newcommand{\Art}{\operatorname{Art}}
\newcommand{\id}{{\rm id}}
\newcommand{\tr}{{\rm tr}}
\newcommand{\ab}{{\rm ab}}
\newcommand{\ur}{{\rm ur}}
\newcommand{\ram}{{\rm ram}}
\newcommand{\cyc}{{\rm cyc}}
\newcommand{\LT}{{\rm LT}}
\newcommand{\isom}{\stackrel{\cong}{\lra}}
\newcommand{\chara}{\operatorname{char}}
\newcommand{\Frac}{\operatorname{Frac}}
\newcommand{\mspec}{\operatorname{m-Spec}}
\newcommand{\Hom}{\operatorname{Hom}}
\newcommand{\res}{\operatorname{res}}
\newcommand{\End}{\operatorname{End}}
\newcommand{\Aut}{\operatorname{Aut}}
\def\bmod{\operatorname{mod}}
\newcommand{\notmid}{\!\not{\mid}\,}
\newcommand{\MU}{\mbox{\boldmath{$\mu$}}}
\newcommand{\sMU}{\mbox{\small\boldmath{$\mu$}}}
\newcommand{\benu}{\begin{enumerate}\renewcommand{\labelenumi}{(\roman{enumi})}\renewcommand{\itemsep}{0pt}}
\newcommand{\eenu}{\end{enumerate}}
\def\hat{\widehat}
\def\invlim#1{\lim\limits_{\substack{\longleftarrow\\#1}}}

\makeatletter
 \renewcommand{\theequation}{%
   \thesubsection.\arabic{equation}}
  \@addtoreset{equation}{subsection}
\makeatother

\title{Local Class Field Theory via Lubin-Tate Theory}
\author{Teruyoshi Yoshida}\thanks{Harvard University, Department of Mathematics, 1 Oxford Street, Cambridge, MA 02138, USA}
\email{yoshida\char`\@math.harvard.edu}

\maketitle

\begin{abstract}
We give a self-contained exposition of local class field theory, via Lubin-Tate theory and the Hasse-Arf theorem, refining the arguments of Iwasawa \cite{Iw}.
\end{abstract}

\def\abstractname{R\'esum\'e}
\begin{abstract}
Nous pr\'esentons une d\'emonstration compl\`ete de la th\'eorie du corps de classes locale via la th\'eorie de Lubin-Tate et le th\'eor\`eme de Hasse-Arf, en raffinant des arguments d'Iwasawa \cite{Iw}.
\end{abstract}

\section{Introduction}

We prove local class field theory via Lubin-Tate theory and the Hasse-Arf theorem. The only prerequisites are Galois theory (including cyclotomic extensions, finite fields and infinite extensions) and some basic commutative algebra summarized in Appendix I. We have tried to make the paper self-contained, to the extent of repeating proofs of standard results on local fields and avoiding topological arguments using compactness. Our argument is close to Iwasawa \cite{Iw}, but the main innovation here is to use the relative Lubin-Tate groups of de Shalit \cite{deSh} to prove the base change property (Theorem \ref{basechange}) directly, without proving the local Kronecker-Weber theorem first.

\begin{ithm} \label{LCFT}
\begin{enumerate}\renewcommand{\labelenumi}{(\roman{enumi})}
\item For any local field $K$, there is a unique homomorphism $\Art_K:K^\times\lra \Gal(K^\ab/K)$, characterized by the two properties:
\begin{enumerate}
\item If $\pi$ is a uniformizer of $K$, then $\Art_K(\pi)|_{K^\ur}=\Frob_K$.
\item If $K'/K$ is a finite abelian extension, then $\Art_K(N_{K'/K}(K^{\prime \times}))|_{K'}=\id$.
\end{enumerate}
Moreover, $\Art_K$ is an isomorphism onto $W_K^\ab:=\{\sigma \mid \sigma|_{K^\ur}\in \Frob_K^\Z\}\subset \Gal(K^\ab/K) $.
\item If $K'/K$ is finite separable, then $\Art_{K'}(x)|_{K^\ab}=\Art_K(N_{K'/K}(x))$ for all $x\in K^{\prime \times}$, and $\Art_K$ induces an isomorphism $K^\times/N_{K'/K}(K^{\prime \times})\isom \Gal((K'\cap K^\ab)/K)$.
\end{enumerate}
\end{ithm}

\paragraph{{\bf Notation}}
The cardinality of a finite set $X$ is denoted by $|X|$. A {\it ring} means a commutative ring with a unit, unless stated otherwise. For a ring $A$, we write $A^\times$ for its group of units. For a field $F$, we usually (implicitly) fix its algebraic closure $\ol{F}$ and separable closure $F^{\rm sep}$, and regard any algebraic (resp.\ separable) extension of $F$ as a subfield of $\ol{F}$ (resp.\ $F^{\rm sep}$). For a finite separable extension $F'/F$, we denote the norm map by $N_{F'/F}:F^{\prime \times}\ra F^\times$. We denote the maximal abelian extension of $F$ in $\ol{F}$ by $F^\ab$. 

For a positive integer $n$ not divisible by $\chara F$, the splitting field of $X^n-1$ over $F$ is denoted by $F(\MU_n)$ ({\it cyclotomic extension}), which is an abelian extension such that its Galois group naturally injects into $(\Z/(n))^\times$. We denote the set of roots of $X^n-1$ by $\MU_n$. For $x\in F^\times$, we write $\langle x \rangle$ for the subgroup $x^\Z:=\{x^a\mid a\in \Z\}$ of $F^\times$ generated by $x$. We denote a finite field consisting of $q$ elements by $\F_q$. For each $n\geq 1$, we have $\F_{q^n}=\F_q(\MU_{q^n-1})$ in $\ol{\F}_q$. The Galois group $\Gal(\ol{\F}_q/\F_q)$ is isomorphic to $\widehat{\Z}:=\invlim{}\Z/(n)$, the {\it profinite completion} of $\Z$, by sending the $q$-th power {\it Frobenius map} $x\mapsto x^q$ to 1.

\section{Local fields and complete extensions} \label{s2}

\subsection{Complete discrete valuation fields (see Appendix I)} \label{cdvf}

Let $K$ be the fraction field of a CDVR $\mO:=\mO_K$ (the {\em ring of integers} of $K$) with maximal ideal $\fkp:=\fkp_K$, such that its {\em residue field} $k:=\mO/\fkp$ is a perfect field. A generator of $\fkp$ is called a {\em uniformizer} of $K$. We denote its {\em valuation} by $v=v_K:K^\times\ra \Z$. If $K'/K$ is a finite separable extension, then $K'$ is the fraction field of a CDVR $\CO_{K'}$, namely the integral closure of $\mO$ in $K'$, and the residue field $k'$ of $K'$ is a finite extension of $k$. The {\em ramification index} $e=e(K'/K)$ and the {\em residue degree} $f=f(K'/K)$ of $K'/K$ are defined by $\fkp\CO_{K'}=\fkp_{K'}^e$ and $[k':k]=f$. Then $[K':K]=ef$, and $v_{K'}|_{K^\times}=e\cdot v_K$ by definition.
If $K''/K'$ is another finite separable extension, clearly $e(K''/K)=e(K''/K')e(K'/K)$ and $f(K''/K)=f(K''/K')f(K'/K)$. We say $K'/K$ is {\em unramified} if $e=1$, and {\em totally ramified} if $f=1$. By the multiplicativity of $e$ and $f$, any intermediate extension of an unramified (resp.\ a totally ramified) extension is again unramified (resp.\ totally ramified). Now for any finite separable $K'/K$, if $F$ is the Galois closure of $K'$, then as the action of $\Gal(F/K)$ preserves $\mO_F$ and hence also $v_F$, we have $v_K(N_{K'/K}(x))=\frac{1}{e(F/K)}v_F(N_{K'/K}(x))=\frac{[K':K]}{e(F/K)}v_F(x)=\frac{[K':K]}{e(K'/K)}v_{K'}(x)=f(K'/K)v_{K'}(x)$ for all $x\in K^{\prime \times}$, i.e.\ we have $v_K\circ N_{K'/K}=f\cdot v_{K'}$.

For any separable extension $E/K$ (not necessarily finite) in $K^{\rm sep}$, the {\em ring of integers} $\mO_E$ of $E$ is defined as the integral closure of $\mO$ in $E$. If $E=\bigcup_{K'}K'$, where $K'/K$ are finite separable, then $\mO_E=\bigcup_{K'}\mO_{K'}$. As $\fkp_{K'}\subset \fkp_{K''}$ whenever $K'\subset K''$, we have an ideal $\fkp_E:=\bigcup_{K'}\fkp_{K'}$ of $\mO_E$, and $\mO_E^\times=\bigcup_{K'} \mO_{K'}^\times=\mO_E\setminus \fkp_E$. Therefore $\mO_E$ is a local ring with the maximal ideal $\fkp_E$, and $E=\bigcup K'=\Frac(\mO_E)$. 

\bde
We call a separable extension $E/K$ {\em unramified} (resp.\ {\em totally ramified}) if it is a union of unramified (resp.\ totally ramified) finite extensions of $K$. We say $E/K$ is {\em finitely ramified} if $E$ is a finite extension of an unramified extension of $K$.
\ede

\blem \label{infinext}
Let $E\subset K^{\rm sep}$ be finitely ramified over $K$. 
\benu
\item The ring of integers $\mO_E$ is a DVR.
\item If $E'/E$ is finite separable, then $E'\hat{E}=\hat{E}'$ and $\hat{E}\cap E'=E$. 
\item $\hat{E}\cap K^{\rm sep}=E$. (Hence $\hat{E}=\hat{E}'\Lra E=E'$ for $E,E'/K$ finitely ramified.)
\eenu
\elem

\bprf
(i): If $E/K$ is unramified, then $\fkp_{K'}=\fkp\mO_{K'}$ for all finite intermediate $K'/K$, therefore $\fkp_E=\fkp\mO_E$ and $\mO_E$ is a DVR. If $E'$ is finite over $E$, then $\mO_{E'}$ is the integral closure of the DVR $\mO_E$ in $E'$, hence a DVR. (ii) follows from Proposition \ref{dvrlemma}(ii). (ii) implies (iii). 
\eprf

\subsection{Local fields and their complete extensions} \label{localfield}

In the rest of the article, we fix a prime $p$, and let $K$ denote a {\em local field}, i.e.\ a complete discrete valuation field whose residue field $k$ is a finite field $\F_q$ of characteristic $p$. Then $\chara K=0$ or $p$, and if $\chara K=0$, then $K$ is a finite extension of the $p$-adic field $\Q_p$. Finite unramified extensions of local fields are classified using the following lemma (see Appendix I for its proof): 

\begin{lem} \label{Hensel}
{\em (Hensel's lemma)} Let $n\geq 1$ with $(p,n)=1$. Then $\MU_n\subset k \iff \MU_n\subset K$. 
\end{lem}

For $n\geq 1$, let $K_n:=K(\MU_{q^n-1})$ and $k_n$ be its residue field. Then $K_n/K$ is unramified (Proposition \ref{cyclounram}), and $\F_{q^n}\subset k_n$ by the above lemma. As $\Gal(K_n/K)\cong \Gal(k_n/\F_q)$ shows that an element of $\Gal(k_n/\F_q)$ is determined by its action on $\MU_{q^n-1}$, we have $k_n=\F_{q^n}$ and $[K_n:K]=n$. Conversely, if $K'/K$ is unramified of degree $n$, then the residue field of $K'$ is $\F_{q^n}$, hence $\MU_{q^n-1}\subset K'$ by the above lemma, and we see $K'=K_n$ by comparing the degrees. As $K_n\subset K_{n'}$ for $n\mid n'$, the union $K^\ur:=\bigcup_{n\geq 1}K_n$ is an infinite Galois extension of $K$ (the {\em maximal unramified extension} of $K$), and by the above isomorphism:
\[\Gal(K^\ur/K)\isom \invlim{}\Gal(\F_{q^n}/\F_q)\cong \Gal(\ol{\F}_q/\F_q)\isom \widehat{\Z}.\]
The {\em arithmetic Frobenius} $\varphi\in \Gal(K^\ur/K)$ is defined as the element which reduces $\bmod \fkp$ to the $q$-th power Frobenius map of $\ol{\F}_q$, and its inverse is denoted by $\Frob_K:=\varphi^{-1}$ ({\em geometric Frobenius}). Unramified extensions of $K$ are none other than subfields of $K^\ur$, hence always abelian over $K$. If $E'/K$ is a separable extension, then $E'$ and $E:=E'\cap K^\ur$ have the same residue fields. When $E'/K$ is Galois, we define its {\em Weil group} by $W(E'/K):=\{\sigma \in \Gal(E'/K)\mid \sigma|_E\in \Frob_K^\Z\}$, which is an extension of $W(E/K)$ (a quotient group of $\Z$) by $\Gal(E'/E)$. If $E/K$ is finite, then $W(E'/K)=\Gal(E'/K)$.

\bde
We call the completion $L=\widehat{E}$ of a finitely ramified (\S\ref{cdvf}) extension $E$ of $K$ a {\em complete extension} of $K$ (if $E/K$ is finite, then $L=E$). Then $\mO_L=\hat{\mO}_E$ is a CDVR with the maximal ideal $\fkp_L=\fkp_E\mO_L$. The complete extensions correspond bijectively to finitely ramified extensions $E/K$ by Lemma \ref{infinext}(iii). When $E/K$ is unramified, we call $L=\widehat{E}$ a {\em complete unramified extension} of $K$. 
\ede

The $\hat{K}:=\hat{K}^\ur$ is a complete unramified extension of $K$, and we write $\hat{\mO}:=\mO_{\hat{K}},\ \hat{\fkp}:=\fkp_{\hat{K}}$. We consider every complete unramified extension $L/K$ as a subfield of $\hat{K}$, in which case $\fkp_L=\fkp\mO_L$ and a uniformizer of $L$ is also a uniformizer of $\hat{K}$. Let $L'=\widehat{E}'$ be a complete extension of $K$, and set $E:=E'\cap K^\ur$. Then $L=\widehat{E}$ is a complete unramified extension of $K$, and $L', E', E, L$ all have the same residue fields, i.e.\ $L'/L$ is totally ramified. We consider every complete extension $L'/K$ as a subfield of $\hat{K}^{\rm sep}$ via $L\subset \hat{K}$.

\bde \label{galweilgp}
Let $L'$ be a totally ramified extension of a complete unramified extension $L/K$. When $L'/L$ is finite, we say $L'$ is {\em Galois over} $K$ if for all $i\in \Z$, the $\varphi^i=\Frob_K^{-i}\in \Aut(L/K)$ extends to $[L':L]$ distinct elements of $\Aut(L'/K)$. In general, we say $L'$ is {\em Galois over} $K$ if it is a union of finite extensions of $L$ which are Galois over $K$. In this case we define the {\em Weil group} of $L'/K$ by $W(L'/K):=\{\sigma \in \Aut(L'/K)\mid \sigma|_L\in \Frob_K^\Z\}$, which is an extension of $W(L/K)$ (a quotient group of $\Z$) by $\Gal(L'/L)$. When $L=\hat{K}$, define $v=v_K:W(L'/K)\ra \Z$ by $\sigma|_L=\Frob_K^{v(\sigma)}$. 
\ede

This terminology coincides with the usual one when $L/K$ is finite. When $L'=\hat{E}'$ for finitely ramified Galois $E'/K$, then every $\sigma\in \Gal(E'/K)$ induces $\mO$-automorphisms of $\mO_{E'}$ and $\mO_{E'}/\fkp_{E'}^m$ for all $m\geq 1$, hence of $\mO_{L'}=\hat{\mO}_{E'}$. Therefore it extends to a $K$-automorphism of $L'$, and we have a canonical injection $\Gal(E'/K)\ra \Aut(L'/K)$. Therefore, as a totally ramified extension of $L=\hat{E}$ for $E=E'\cap K^\ur$, we see that $L'$ is Galois over $K$ (because $[L':L]=[E':E]$ by Lemma \ref{infinext}(ii)), and canonically $W(E'/K)\cong W(L'/K)$. By passing to the limit, this last isomorphism extends to the case where $L'=E'L$ with a general Galois extension $E'/K$.

\section{Formal groups and Lubin-Tate groups}

\subsection{Formal groups}

Let $A$ be a ring, not the zero ring. In the formal power series ring of one variable $A[[X]]:=\invlim{m}A[X]/(X^m)$ over $A$, the ideal $(X)\subset A[[X]]$, consisting of all the elements with constant term equal to 0, is a monoid under the composition $f\circ g:=f(g(X))$ with $X$ as the identity. For $f\in (X)$, there exists an $f^{-1}$ satisfying $f\circ f^{-1}=f^{-1}\circ f=X$ if and only if the coefficient of $X$ in $f$ belongs to $A^\times$. Also, we use similar notation for $f\in (X)\subset A[[X]]$ and a power series of several variables $F\in A[[X_1,\ldots,X_n]]$ with no constant term:
\[f\circ F:=f(F(X_1,\ldots,X_n)),\ \ \ F\circ f:=F(f(X_1),\ldots,f(X_n))\in A[[X_1,\ldots,X_n]].\]

\begin{defn}
A {\em formal group over $A$} is a formal power series of two variables $F(X,Y)\in A[[X,Y]]$ which satisfies the following: 
\benu
\item $F(X,Y)\equiv X+Y\ (\bmod \deg 2)$,
\item $F(F(X,Y),Z)=F(X,F(Y,Z))$,
\item $F(X,Y)=F(Y,X)$.
\eenu
\end{defn}

Precisely speaking, these are {\em commutative formal groups of dimension 1}. The basic examples are the {\em additive group} $\widehat{\G}_a(X,Y):=X+Y$ and the {\em multiplicative group} $\widehat{\G}_m(X,Y):=X+Y+XY$.

Let $F$ be a formal group over a ring $A$. If we let $f(X):=F(X,0)$, we have $f(X)\equiv X\ (\bmod \deg 2)$ by (i), hence $f^{-1}$ exists. By (ii), we have $f\circ f=f$, hence we get $f(X)=X$ by composing with $f^{-1}$. Similarly we have $F(0,Y)=Y$, hence $F$ does not have a term containing only $X$ or $Y$, apart from the linear terms $X+Y$. Therefore we can solve $F(X,Y)=0$ with respect to $Y$ and get a unique $i_F(X)\in A[[X]]$ satisfying $F(X,i_F(X))=0$. If we define the addition $+_F$ on the ideal $(X)\subset A[[X]]$ by
\[f+_Fg:=F(f(X),g(X)),\]
then $(X)$ becomes an abelian group with 0 as the identity and $i_F\circ f$ as the inverse of $f$.

\begin{defn}
Let $F,G$ be formal groups over $A$. A power series $f(X)\in (X)\subset A[[X]]$ is called a {\em homomorphism} from $F$ to $G$ if it satisfies 
\[f\circ F=G\circ f,\ \ \text{i.e.}\ \ f(F(X,Y))=G(f(X),f(Y)),\]
and we write $f:F\ra G$. Two homomorphisms compose via the composition of power series, with $f(X)=X$ as the identity $\id:F\ra F$. If $f^{-1}$ exists, it defines $f^{-1}:G\ra F$ and $f\circ f^{-1}=f^{-1}\circ f=\id$. In this case $f$ is called an {\em isomorphism} and we write $f:F\isom G$. 
\end{defn}

The set $\Hom_A(F,G)$ of all homomorphisms from $F$ to $G$ is an abelian group under $+_G$. Moreover, $\End_A(F):=\Hom_A(F,F)$ is a (not necessarily commutative) ring with $+_F$ as the addition and $\circ$ as the multiplication.


\subsection{Lubin-Tate groups}

We return to the notation of \S\ref{localfield}, i.e.\ $K$ is a local field with the ring of integers $\mO$ and its maximal ideal $\fkp$, and $\mO/\fkp\cong \F_q$ where $q$ is a power of $p$. Let $L$ be a {\em complete unramified extension} of $K$ (\S\ref{localfield}). As $\fkp_L=\fkp\mO_L$, we write $\bmod \fkp$ for $\bmod \fkp\mO_L$. Let $\varphi$ be the arithmetic Frobenius, extended to a $K$-automorphism of $L$. For $\alpha\in L$ and $i\in \Z$, we write $\alpha^{\varphi^i}:=\varphi^i(\alpha)$. For a power series $F$ over $\mO_L$, we define $F^{\varphi^i}$ by applying $\varphi^i$ to all coefficients of $F$. If $F$ is a formal group over $\mO_L$, so is $F^{\varphi^i}$.

\bde
For uniformizers $\pi,\pi'$of $L$, set $\Theta^L_{\pi,\pi'}:=\{\theta\in \mO_L\mid \theta^\varphi/\theta = \pi'/\pi \}$. It is an additive group. If $\theta \in \Theta^L_{\pi,\pi'}$ and $\theta'\in \Theta^L_{\pi',\pi''}$, then $\theta\theta'\in \Theta^L_{\pi,\pi''}$. We have $\mO\subset \Theta^L_{\pi,\pi}$ (actually we will see $\mO=\Theta_{\pi,\pi}^L$ by Lemma \ref{unramnorm}(i)). 
\ede

\blem \label{fundlem4}
Let $\pi$ be a uniformizer of $L$, and let $f\in \mO_L[[X]]$ satisfy:
\begin{equation} \label{fpi}
f(X)\equiv \pi X\ (\bmod \deg 2),\ \ \ \ f(X)\equiv X^q\ (\bmod \fkp).
\end{equation}
Let $\pi',f'$ be another such pair. Assume that $\theta_1,\ldots,\theta_t\in \Theta_{\pi,\pi'}^L$. Then there is a unique $F\in \mO_L[[X_1,\ldots,X_t]]$ satisfying the following: 
\[F\equiv \theta_1X_1+\cdots +\theta_tX_t\ (\bmod \deg 2),\ \ \ \ f'\circ F=F^\varphi\circ f.\]
\elem

\bprf
It suffices to show that for each $m\geq 1$, there is a unique polynomial $F_m$ of degree $\leq m$ that satisfies the conditions $\bmod \deg (m+1)$. The case $m=1$ is assumed, and suppose we have $F_m$, and let $G_{m+1}:=f'\circ F_m-F_m^\varphi\circ f$. Then as $G_{m+1}\equiv F_m^q-F_m^\varphi(X_1^q,\ldots,X_n^q)\equiv 0\ (\bmod \fkp)$, its coefficients are divisible by $\pi'$. Now we show that a homogeneous polynomial $H_{m+1}:=F_{m+1}-F_m$ of degree $m+1$ is uniquely determined. We need $f'\circ F_{m+1}-F_{m+1}^\varphi\circ f\equiv G_{m+1}+(f'\circ H_{m+1}-H_{m+1}^\varphi\circ f)\equiv G_{m+1}+(\pi'H_{m+1}-\pi^{m+1}H_{m+1}^\varphi)\ (\bmod \deg (m+2))$ to vanish. For any monomial of degree $m+1$, if we let $\pi'\beta$ be its coefficient in $G_{m+1}$, and $\alpha$ its coefficient in $H_{m+1}$, then $\pi'\beta+\pi'\alpha-\pi^{m+1}\alpha^\varphi=0$, hence $\alpha=-\beta-\sum_{i=1}^\infty(\pi^{m+1}/\pi')^{1+\varphi+\cdots \varphi^{i-1}}\beta^{\varphi^i}$.
\eprf

\bpr \label{LTexists}
Let $f,f'\in \mO_L[[X]]$ be as above, with linear coefficients $\pi,\pi'$ respectively. 
\benu
\item There exists a unique formal group $F_f$ over $\mO_L$ such that $f\in \Hom_{\mO_L}(F_f,F_f^\varphi)$. (We call $F_f$ the {\em Lubin-Tate group} associated to $f$.)
\item There is a unique map $[\cdot]_{f,f'}:\Theta^L_{\pi,\pi'}\ra (X)\subset \mO_L[[X]]$ such that:
\[ [\theta]_{f,f'}(X)\equiv \theta X\ (\bmod \deg 2),\ \ \ \ f'\circ [\theta]_{f,f'}=[\theta]_{f,f'}^\varphi\circ f.\]
It satisfies $[\theta]_{f,f'}+_{F_{f'}}[\theta']_{f,f'}=[\theta+\theta']_{f,f'},\ \ [\theta']_{f',f''}\circ [\theta]_{f,f'}=[\theta\theta']_{f,f''}$. 
\item We have $[\theta]_{f,f'}\in \Hom_{\mO_L}(F_f,F_{f'})$ for all $\theta\in \Theta^L_{\pi,\pi'}$. 
\eenu
\epr

\bprf 
(i): Lemma \ref{fundlem4} for $\pi=\pi',\ f=f',\ t=2,\ \theta_1=\theta_2=1$ gives a unique $F_f\in \mO_L[[X,Y]]$ with $F_f\equiv X+Y\ (\bmod \deg 2)$ and $f\circ F_f=F_f^\varphi \circ f$. As $F_f(Y,X)$ enjoys the same property, $F_f(X,Y)=F_f(Y,X)$. Similarly, $F_f(F_f(X,Y),Z)$ and $F_f(X,F_f(Y,Z))$ both satisfy the conditions of the lemma for $t=3$ and $\theta_1=\theta_2=\theta_3=1$, hence are equal. Thus $F_f$ is a formal group and $f\in \Hom_{\mO_L}(F_f,F_f^\varphi)$. (ii): Lemma \ref{fundlem4} for $t=1$ gives $[\theta]_{f,f'}$. 
The properties characterizing $[\theta+\theta']_{f,f'}$ (resp.\ $[\theta\theta']_{f,f''}$) are shared by $[\theta]_{f,f'}+_{F_{f'}}[\theta']_{f,f'}$ (resp.\ $[\theta']_{f',f''}\circ [\theta]_{f,f'}$) because:
\begin{gather*}
f'\circ ([\theta]+_{F_{f'}}[\theta']) = (f'\circ [\theta])+_{F_{f'}^\varphi}(f'\circ [\theta'])=([\theta]^\varphi +_{F_{f'}^\varphi} [\theta']^\varphi)\circ f = ([\theta]+_{F_{f'}}[\theta'])^\varphi\circ f\\
(\text{resp. } f''\circ ([\theta']\circ [\theta]) = [\theta']^\varphi \circ f' \circ [\theta]) = [\theta']^\varphi \circ [\theta]^\varphi \circ f = ([\theta']\circ [\theta])^\varphi \circ f\ ).
\end{gather*}
(iii): For $[\theta]:=[\theta]_{f,f'}$, we have $[\theta]\circ F_f=F_{f'}\circ [\theta]$, because the equalities:
\begin{gather*}
f'\circ ([\theta]\circ F_f)=[\theta]^\varphi\circ f\circ F_f=([\theta]^\varphi\circ F_f^\varphi)\circ f=([\theta]\circ F_f)^\varphi\circ f,\\
f'\circ (F_{f'}\circ [\theta])=F_{f'}^\varphi\circ f'\circ [\theta]=(F_{f'}^\varphi\circ [\theta]^\varphi)\circ f=(F_{f'}\circ [\theta])^\varphi\circ f,
\end{gather*}
show that both sides satisfy the conditions of Lemma \ref{fundlem4} for $\pi=\pi',\ t=2,\ \theta_1=\theta_2=\theta$. 
\eprf

\bexa \label{cycloaslt}
If $K=\Q_p,\ \pi=p$ and $f=(1+X)^p-1$, then $F_f=\hat{\G}_m=X+Y+XY$.
\eexa

\bcor \label{formalOmodisom}
\benu
\item The map $[\cdot]_f:=[\cdot]_{f,f}:\mO\lra \End_{\mO_L}(F_f)$ is an injective ring homomorphism. (Hence $(F_f,[\cdot]_f)$ is a {\em formal $\mO$-module}.)
\item If $\theta\in \Theta^{L,\times}_{\pi,\pi'}:=\Theta^L_{\pi,\pi'}\cap \mO_L^\times$, then $[\theta]_{f,f'}$ is an isomorphism with the inverse $[\theta^{-1}]_{f',f}$.
\eenu
\ecor

\bexa \label{fpimult}
We have $\pi\in \Theta^L_{\pi,\pi^\varphi}$, and $[\pi]_{f,f^\varphi}=f:F_f\ra F_f^\varphi$ for $f$ satisfying (\ref{fpi}), by uniqueness. (Also note that $F_f^\varphi=F_{f^\varphi}$ and $[\theta]_{f,f'}^\varphi=[\theta^\varphi]_{f^\varphi,f^{\prime \varphi}}$ by uniquness.)
\eexa

\bde \label{fmpim} 
Generalizing Example \ref{fpimult}, define $f_m:=f^{\varphi^{m-1}}\circ\cdots\circ f^\varphi\circ f\in \mO_L[[X]]$ for $m\geq 1$, and set $f_0(X):=X$. Then, by Example \ref{fpimult} and Proposition \ref{LTexists}(ii):
\[ f_m=[\pi^{\varphi^{m-1}}]_{f^{\varphi^{m-1}},f^{\varphi^m}}\circ\cdots \circ [\pi^\varphi]_{f^\varphi,f^{\varphi^2}}\circ [\pi]_{f,f^{\varphi}}=[\pi_m]_{f,f^{\varphi^m}}\ \ (\forall m\geq 0), \]
where we define $\pi_m\in \mO_L$ by $\pi_m:=\prod_{t=0}^{m-1}\pi^{\varphi^t}$ and $\pi_0:=1$.
\ede

\section{Lubin-Tate extensions and Artin maps}

\subsection{Lubin-Tate extensions}
Here we fix a complete unramified extension $L$ of $K$. 

\bde \label{fmdef}
Let $f\in \mO_L[X]$ be a monic polynomial satisfying (\ref{fpi}) for a uniformizer $\pi$ of $L$. For $m\geq 1$, let $L_f^m$ be the splitting field of $f_m\in \mO_L[X]$ (Definition \ref{fmpim}) over $L$, and let $\MU_{f,m}:=\{\alpha\in L_f^m\mid f_m(\alpha)=0\}.$
\ede

\bexa \label{cycloaslt2}
In Example \ref{cycloaslt}, we have $f_m(X)=[p^m]_f(X)=(1+X)^{p^m}-1,\ \MU_{f,m}=\{\zeta-1\mid \zeta\in \MU_{p^m}\}$ and $L_f^m=L(\MU_{p^m})$ for all $m\geq 0$.
\eexa

\blem \label{mupmtors}
Let $m\geq 1$ and $f\in \mO_L[X]$ as above, and set $L':=L_f^m$ and $[\cdot]:=[\cdot]_f$.
\benu
\item The extension $L'/L$ is separable and $\MU_{f,m}\subset \fkp_{L'}$. (In particular, we can substitute the elements of $\MU_{f,m}$ into power series over $\mO_L$ (see Appendix I).)
\item For $x\in K^\times$ with $v(x)=m$ and $\alpha \in \fkp_{L^{\rm sep}}$:
\[ \alpha\in \MU_{f,m} \iff [x](\alpha)=0\iff [a](\alpha)=0\ (\forall a\in \fkp^m).\] 
\eenu
\elem

\bprf (i): The separability of $L'/L$ is automatic when $\chara K=0$, and in general it follows from Proposition \ref{fmsep} in the Appendix II (which in turn follows from Proposition \ref{LTextension}(i) when $\chara K=0$). Now $\MU_{f,m}\subset \mO_{L'}$ as $f_m$ is a monic in $\mO_L[X]$. If $\alpha\in \mO_{L'}^\times$, then $f_m(\alpha)$, being $\equiv \alpha^{q^m} (\bmod \fkp_{L'})$, will also be in $\mO_{L'}^\times$. Thus $\MU_{f,m}\subset \fkp_{L'}$. (ii): By Definition \ref{fmpim}, we have $[x]=[x/\pi_m]_{f^{\varphi^m},f}\circ f_m$. As $[x/\pi_m]$ is invertible, we see the first equivalence. The second one follows by $\fkp^m=(x)$.
\eprf

\bpr \label{LTextension}
Let $m\geq 1$ and $f\in \mO_L[X]$ as above, with the linear coefficient $\pi$.
\benu
\item The set $\MU_{f,m}$ is an $\mO$-module by $+_{F_f}$ and $[\cdot]_f$. For any $\alpha\in \MU_{f,m}^\times:= \MU_{f,m}\setminus \MU_{f,m-1}$, the following is an isomorphism of $\mO$-modules: 
\[\mO/\fkp^m\ni a\bmod \fkp^m \longmapsto [a]_f(\alpha)\in \MU_{f,m}.\]
\item If $\alpha\in \MU_{f,m}^\times$, then $L_f^m=L(\alpha),\ N_{L_f^m/L}(-\alpha)=\pi^{\varphi^{m-1}}$ and $\alpha$ is a uniformizer of $L_f^m$. The $L_f^m/L$ is totally ramified Galois extension of degree $|\MU_{f,m}^\times|=q^{m-1}(q-1)$. 
\item We have canonical isomorphisms of abelian groups: 
\begin{align*}
\rho_{f,m}: \Gal(L_f^m/L) \isom \Aut_\mO(\MU_{f,m}) &\isom (\mO/\fkp^m)^\times.\\
(\alpha\mapsto [u]_f(\alpha),\ \forall \alpha\in \MU_{f,m}) &\longmapsto u\bmod \fkp^m
\end{align*}
\eenu
\epr

\bprf
We write $+_f:=+_{F_f}$ and $L':=L_f^m$. (i): Lemma \ref{mupmtors}(ii) shows that $\MU_{f,m}$ is an $\mO$-module by $+_f,[\cdot]$, killed by $\fkp^m$. The stated $\mO$-homomorphism is injective as $[a](\alpha)\neq 0$ for some $a\in \fkp^{m-1}$ by Lemma \ref{mupmtors}(ii), hence surjective as $|\mO/\fkp^m|=q^m=\deg f_m\geq |\MU_{f,m}|$. (Thus $|\MU_{f,m}|=q^m$ and hence $\MU_{f,m}^\times$ is the set of all roots of $f_m/f_{m-1}$.) (ii): We have $\MU_{f,m}\subset L(\alpha)$ by (i), hence $L'=L(\alpha)$ and $L'/L$ is Galois. Now the constant term of $f_m/f_{m-1}$ reads $\pi^{\varphi^{m-1}}=\prod_{\alpha\in \sMU_{f,m}^\times}(-\alpha)$, and taking the $v_{L'}$ of both sides shows $e(L'/L)=\sum v_{L'}(-\alpha)\geq |\MU_{f,m}^\times|$ by Lemma \ref{mupmtors}(i). But $|\MU_{f,m}^\times|=\deg(f_m/f_{m-1})\geq [L':L]\geq e(L'/L)$, hence all are equalities and $f_m/f_{m-1}$ is irreducible. (iii): As $+_f,[\cdot]$ have coefficients in $\mO_L$, for all $\sigma\in \Gal(L'/L)$, we have $\sigma(\alpha +_f \alpha')=\sigma(\alpha)+_f\sigma(\alpha')$ and $\sigma([a](\alpha))=[a](\sigma(\alpha))$, i.e.\ $\Gal(L'/L)$ acts on $\MU_{f,m}$ by $\mO$-homomorphisms. Hence we have a group homomorphism $\rho_{f,m}:\Gal(L'/L)\lra \Aut_\mO(\MU_{f,m})$. This is injective as $L'=L_f^m$, and $\Aut_\mO(\MU_{f,m})\cong (\mO/\fkp^m)^\times$ by (i). It is surjective as $|\Gal(L'/L)|=[L':L]=|(\mO/\fkp^m)^\times|$ by (ii).
\eprf

\subsection{Artin map}

In this subsection we use the notation $(\ )^{(i)}:=(\ )^{\varphi^i}$ and $\MU_{f,m}^{(i)}:=\MU_{f^{(i)},m}$ for all $i\in \Z$. We extend Definition \ref{fmpim} to define $\pi_j\in L^\times$ for all $j\in \Z$ by requiring $\pi_{j+j'}=\pi_{j'}^{(j)}\pi_j$ for all $j,j'\in \Z$, i.e.\ $\pi_j:=(\pi_{-j}^{-1})^{(j)}$ for $j<0$. Then $v_L(\pi_j)=j$ for all $j\in \Z$. 

\blem \label{thetapi}
If $\theta\in \Theta_{\pi,\pi'}^L$, then $\theta^{(j)}/\theta=\pi'_j/\pi_j$ for all $j\in \Z$. Also, $\pi_j\in \Theta_{\pi,\pi^{(j)}}^L$.
\elem

\bprf
Using $\pi'_{j+1}/\pi_{j+1}=(\pi'_j/\pi_j)(\pi'/\pi)^{(j)}=(\pi'_j/\pi_j)(\theta^\varphi/\theta)^{(j)}=(\pi'_j/\pi_j)(\theta^{(j+1)}/\theta^{(j)})$, argue by induction in both directions. Take $\pi'=\pi^\varphi$ and $\theta=\pi$ for the second claim.
\eprf

\blem \label{isomthensamefield}
Let $f,f'\in \mO_L[X]$ be as above with linear coefficients $\pi,\pi'$, respectively. If $\theta\in \Theta_{\pi,\pi'}^{L,\times}$ (see Corollary \ref{formalOmodisom}(ii)), then for all $m\geq 1$, it gives an isomorphism $[\theta]=[\theta]_{f,f'}:\MU_{f,m}\ra \MU_{f',m}$ of $\mO$-modules, and $L_f^m=L_{f'}^m$. 
\elem

\bprf
The $[\theta]$ maps $\MU_{f,m}$ to $\MU_{f',m}$ because $f'_m\circ [\theta]=[\theta]^{(m)}\circ f_m$. It is an $\mO$-homomorphism by Proposition \ref{LTexists}(ii),(iii), and is an isomorphism as $[\theta^{-1}]$ gives its inverse. As $[\theta],[\theta^{-1}]\in \mO_L[[X]]$, we have $\MU_{f',m}=[\theta](\MU_{f,m})\subset L_f^m$ and $\MU_{f,m}\subset L_{f'}^m$, thus $L_f^m=L_{f'}^m$. 
\eprf

\bpr \label{invartmap}
Let $m\geq 1$ and $f\in \mO_L[X]$ as above, with the linear coefficient $\pi$.
\benu
\item The $L_f^m$ is Galois over $K$, and the following map is bijective for any $\alpha\in \MU_{f,m}^\times$:
\[ K^\times / (1+\fkp^m)\ni x\bmod 1+\fkp^m \longmapsto [x\pi_j]_{f,f^{(j)}}(\alpha)\in \coprod_{j\in \Z}\MU_{f,m}^{(j),\times}\ \ (v(x)=-j).\]
\item Let $L=\hat{K}$. The $\rho_{f,m}$ of Proposition \ref{LTextension}(iii) extend to isomorphisms: 
\begin{align*}
\rho_{f,m}: W(\hat{K}_f^m/K) &\isom K^\times / (1+\fkp^m).\\
(\varphi^j \text{ on } \hat{K},\ \ \alpha\mapsto [x\pi_j](\alpha),\ \forall \alpha\in \MU_{f,m}) &\longmapsto x\bmod 1+\fkp^m\ \ (v(x)=-j)
\end{align*}
Setting $\hat{K}^\LT_f:=\bigcup_{m\geq 1} \hat{K}_f^m$, we get $\rho_f:W(\hat{K}^\LT_f/K)\isom K^\times$ by passing to the limit. 
\eenu
\epr

\bprf
(i): If $v(x)=-j$, then $x\pi_j\in \Theta^{L,\times}_{\pi,\pi^{(j)}}$ by Lemma \ref{thetapi}, hence $[x\pi_j]:\MU_{f,m}\isom \MU_{f,m}^{(j)}$ by Lemma \ref{isomthensamefield}. As $[x\pi_j]$ is $\mO$-linear, $v^{-1}(-j)/(1+\fkp^m)\ni x\mapsto [x\pi_j](\alpha)\in \MU_{f,m}^{(j),\times}$ is bijective for each $j$. As $L(\alpha)=L_f^m=L_{f^{(j)}}^m$ by Proposition \ref{LTextension}(ii) and Lemma \ref{isomthensamefield}, the $\varphi^j\in \Aut(L/K)$ extends to $L_f^m$ by $\alpha\mapsto \alpha'$ for each $\alpha'\in \MU_{f,m}^{(j),\times}$, hence $L_f^m$ is Galois over $K$. (ii): Let $\sigma\in W(\hat{K}_f^m/K)$ with $\sigma|_{\hat{K}}=\varphi^j$. If $\alpha\in \MU_{f,m}^\times$, then $\sigma(\alpha)\in \MU_{f,m}^{(j),\times}$, hence $\sigma(\alpha)=[x\pi_j](\alpha)$ for a unique $x \bmod 1+\fkp^m$ by (i). This holds for all $\alpha\in \MU_{f,m}$ because $\sigma([a]_f(\alpha))=[a]_f^{(j)}(\sigma(\alpha))=[a]_{f^{(j)}}[x\pi_j](\alpha)=[x\pi_j][a]_f(\alpha)$ for all $a\in \mO$ (this shows the compatibility of $\rho_{f,m}$ for varying $m$). The map $\rho_{f,m}$ is a group homomorphism because if $\tau(\alpha)=[y\pi_{j'}](\alpha)$, then $\sigma\tau(\alpha)=\sigma([y\pi_{j'}](\alpha))=[y\pi_{j'}]^{(j)}[x\pi_j](\alpha)=[y\pi_{j'}^{(j)}\cdot x\pi_j](\alpha)=[xy\cdot \pi_{j+j'}](\alpha)$. It is bijective because it restricts to $\Gal(\hat{K}_f^m/\hat{K})\cong (\mO/\fkp^m)^\times=\mO^\times/(1+\fkp^m)$ by Proposition \ref{LTextension}(iii) and the quotient $W(\hat{K}/K)=\Frob_K^\Z$ is mapped onto $K^\times/\mO^\times\cong \Z$, i.e.\ $v\circ \rho_{f,m}=v$.
\eprf

\bpr \label{thetafunct2}
The map $\psi:\hat{\mO}^\times\ni \theta\longmapsto \theta^\varphi/\theta\in \hat{\mO}^\times$ is surjective. In particular, for any pair of uniformizers $\pi,\pi'$ of $\hat{K}$, we have $\Theta_{\pi,\pi'}^{\hat{K},\times} \neq \emptyset$.
\epr

\bprf
As $\hat{\mO}^\times\cong \invlim{}(\hat{\mO}/\hat{\fkp}^m)^\times= \invlim{}\hat{\mO}^\times/(1+\hat{\fkp}^m)$ and $\psi(1+\hat{\fkp}^m)\subset 1+\hat{\fkp}^m$, it suffices to show for every $u\in \hat{\mO}^\times$ and all $m\geq 1$, there is $\theta_m\in \hat{\mO}^\times$ with $\psi(\theta_m)\equiv u\ (\bmod \fkp^m)$ and $\theta_{m+1}\equiv \theta_m\ (\bmod \fkp^m)$. We get $\theta_1$ because $\ol{\theta}\mapsto \ol{\psi(\theta)}=\ol{\theta}^{q-1}$ is surjective on $(\hat{\mO}/\hat{\fkp})^\times\cong \ol{\F}_q^\times$. Suppose we have $\theta_m$, and let $u/\psi(\theta_m)=1+\alpha\pi^m$ for a uniformizer $\pi$ of $K$. Then there is $\beta\in \hat{\mO}$ with $\beta^\varphi-\beta\equiv \alpha\ (\bmod \fkp)$ because $\ol{\beta}\mapsto \ol{\beta^\varphi}-\ol{\beta}=\ol{\beta}^q-\ol{\beta}$ is surjective on $\hat{\mO}/\hat{\fkp}\cong \ol{\F}_q$, and $\theta_{m+1}:=\theta_m(1+\beta\pi^m)$ will do.
\eprf

\bcor
The $\hat{K}_f^m$ and $\rho_{f,m}$, hence also $\hat{K}^\LT_f$ and $\rho_f$, of Proposition \ref{invartmap}(ii) do not depend on $f$. (We will drop the subscript $f$ and write $\hat{K}^m,\ \rho_m,\ \hat{K}^\LT$ and $\rho$.)
\ecor

\bprf
For $f,f'$ with linear coefficients $\pi,\pi'$, take $\theta\in \Theta_{\pi,\pi'}^{\hat{K},\times}$ and $[\theta]:\MU_{f,m}\isom \MU_{f',m}$ by Proposition \ref{thetafunct2}. Lemma \ref{isomthensamefield} shows $\hat{K}_f^m=\hat{K}_{f'}^m$. If $\sigma(\alpha)=[x\pi_j](\alpha)$ for $\sigma\in W(\hat{K}_f^m/K)$, then $\sigma([\theta](\alpha)) = [\theta]^{(j)}[x\pi_j](\alpha)= [x\pi'_j][\theta](\alpha)$ by Lemma \ref{thetapi}, hence $\rho_{f,m}=\rho_{f',m}$.
\eprf

\bde \label{artmap}
For any $f\in \mO_L[X]$ with $L/K$ finite, set $K^m:=K^\ur L_f^m$. Then $K^m/K$ is finitely ramified, and Galois by Proposition \ref{invartmap}(i). By Lemma \ref{infinext}, the completion of $K^m$ is $\hat{K}L_f^m=\hat{K}^m$ and $K^m=\hat{K}^m\cap K^{\rm sep}$, thus independent of $f$. Setting $K^\LT:=\bigcup_{m\geq 1}K^m=\hat{K}^\LT\cap K^{\rm sep}$, we have 
$W(K^\LT/K)\cong W(\hat{K}^\LT/K)$ by the remark after Definition \ref{galweilgp}. We call a finite extension of $K$ a {\em Lubin-Tate extension} if it is contained in $K^\LT$. We call the inverse of $\rho$ the {\em Artin map} of $K$ and write $\Art_K:K^\times \isom W(K^\LT/K)$. We have $v\circ \Art_K=v$. 
\ede

\section{Galois Groups, Norm Groups and the Base Change} \label{LTbasechange}

\subsection{Galois groups}

Now let $L=K_n/K$ be the {\em finite} unramified extension of degree $n$. 

\bpr \label{normisomfin}
Let $\theta\in \Theta_{\pi,\pi'}^{\hat{K},\times}$ for $\pi,\pi'\in L$. Then $\theta\in \mO_L^\times \iff N_{L/K}(\pi)=N_{L/K}(\pi')$.
\epr

\bprf
Lemma \ref{thetapi} for $j=n$ shows $\theta^{\varphi^n}/\theta=N_{L/K}(\pi')/N_{L/K}(\pi)$, so use Lemma \ref{unramnorm}(i).
\eprf

\blem \label{unramnorm}
\benu
\item For $n\geq 1$, the fixed field of $\varphi^n$ in $\hat{K}$ is $K_n$.
\item If $L=K_n$, then $N=N_{L/K}$ surjects onto $v^{-1}(n\Z)\subset K^\times$.
\eenu
\elem

\bprf
(i): As a set of representatives of $\hat{\mO}/\hat{\fkp}\cong \ol{\F}_q$, we can take $C:=\{0\}\cup \bigcup_{n\geq 1}\MU_{q^n-1}$ by Lemma \ref{Hensel}. Then $\varphi^n$ acts on $C$, and its fixed set is $C_n=\{0\}\cup \MU_{q^n-1}\subset K_n$. Now take a uniformizer $\pi$ of $K$, and consider the $\pi$-adic expansion in $\hat{K}$ with respect to $C$ (see Appendix I). If $x=\sum_{i=v(x)}^\infty a_i\pi^i$ for $a_i\in C$, then $x^{\varphi^n}=\sum_ia_i^{\varphi^n}\pi^i$, hence $x^{\varphi^n}=x\iff a_i\in C_n\ (\forall i)\iff x\in K_n$. (ii): For a uniformizer $\pi$ of $K$, we have $v^{-1}(n\Z)=\mO^\times\times \langle \pi^n \rangle$ and $N(\pi)=\pi^n$, hence it suffices to show that $N:\mO_L^\times \ra \mO^\times$ is surjective. We have $\mO^\times\cong \invlim{}\mO^\times/(1+\fkp^m)$, $\mO_L^\times\cong \invlim{}\mO_L^\times/(1+\fkp_L^m)$, and $N(1+\fkp_L^m)\subset 1+\fkp^m$, because $N(1+\fkp_L^m)\subset (1+\fkp_L^m)\cap \mO=1+\fkp^m$. Therefore it suffices to show that, for every $x\in \mO^\times$ and all $m\geq 1$, there is $u_m\in \mO_L^\times$ satisfying $N(u_m)\equiv x\ (\bmod \fkp^m)$ and $u_{m+1}\equiv u_m\ (\bmod \fkp^m)$. We get $u_1$ by the surjectivity of the norm map $(\mO_L/\fkp_L)^\times\ra (\mO/\fkp)^\times$ induced by $N$. Suppose we have $u_m$, and let $x/N(u_m)=1+\alpha\pi^m$. Then there is $\beta\in \mO_L$ whose trace $\equiv \alpha\ (\bmod \fkp)$ because the trace map $\mO_L/\fkp_L \ra \mO/\fkp$ is surjective, and $u_{m+1}:=u_m(1+\beta\pi^m)$ will do.
\eprf

\bde \label{finLTkxm}
Let $x\in K^\times$ with $v(x)=n>0$. Take a uniformizer $\pi$ of $L=K_n$ with $N_{L/K}(\pi)=x$ by Lemma \ref{unramnorm}(ii), and a monic $f\in \mO_L[X]$ satisfying (\ref{fpi}) for $\pi$. Then for $m\geq 1$, the fields $L_f^m$ depend only on $x$ by Proposition \ref{normisomfin} and Lemma \ref{isomthensamefield}, so we denote them by $K_x^m:=L_f^m$, and set $K_x^\ram:=\bigcup_{m\geq 1}K_x^m$, which are totally ramified over $L$. 
\ede

\bpr \label{artchar}
For $x\in K^\times$ with $v(x)=n>0$, the element $\sigma:=\Art_K(x)\in W(K^\LT/K)$ is characterized by $v(\sigma)=v(x)$ and $\sigma|_{K_x^\ram}=\id$. For all $m\geq 1$, the Artin map induces the isomorphism $K^\times / \bigl((1+\fkp^m)\times \langle x\rangle \bigr) \isom \Gal(K_x^m/K)$. 
\epr

\bprf
The $\sigma$ acts as $\Frob_K^n$ on $L$, and $\pi_n=x$ implies $[x\pi_{-n}]=[1]=\id$ on $\MU_{f,m}$, hence $\sigma$ fixes $K_x^\ram$. This characterizes $\sigma$ because $K^\LT=K^\ur K_x^\ram$. It also shows that $\Art_K$ (or $\rho_m^{-1}$) descends to the claimed map, which is bijective because it restricts to $(\mO/\fkp^m)^\times\cong \Gal(K_x^m/L)$ and induces $K^\times/\bigl(\mO^\times\times \langle x\rangle \bigr)\cong \Gal(L/K)$ on the quotients, as $v\circ \Art_K=v$.
\eprf

\subsection{Coleman operator and norm groups}

As above, let $f\in \CO_L[X]$ be a monic polynomial satisfying (\ref{fpi}) for a uniformizer $\pi$ of $L=K_n$, and set $x:=N_{L/K}(\pi)$. We write $+_f$ for $+_{F_f}$ and $\MU_m$ for $\MU_{f,m}$ (we will not see roots of unity here), so $K_x^m=L(\MU_m)$. 

\begin{lem} \label{picompimg}
Let $g\in \CO_L[[X]]$.
\benu
\item If $g(\alpha)=0$ for all $\alpha\in \MU_1$, then $g=g'\cdot f$ for some $g'\in \CO_L[[X]]$.
\item For $h\in \CO_L[[X]]$ and $m\geq 1$, we have $h\circ f\equiv 0\ (\bmod \fkp^m) \Lra h\equiv 0\ (\bmod \fkp^m)$. 
\item If $g(X+_f\alpha)=g(X)$ for all $\alpha\in \MU_1$, then $g=h\circ f$ for a unique $h\in \CO_L[[X]]$.
\eenu
\end{lem}

\pfbegin
(i): For $\alpha\in \MU_1$, if $g(X)=\sum_{i=0}^\infty a_iX^i,\ g(\alpha)=0$ then if we let $b_i:=\sum_{j=0}^\infty a_{i+j+1}\alpha^j\in \CO_{L'}$ for each $i\geq 0$, then $g(X)=(X-\alpha)\cdot\sum_{i=0}^\infty b_iX^i$ in $\CO_{L'}[[X]]$. As $f$ is separable (by Proposition \ref{fmsep}, or Proposition \ref{LTextension}(i) when $\chara K=0$), repeating this, we get $g(X)=f(X)\cdot g'(X)$, and as $g,f \in \CO_L[[X]]$, also $g'$ has coefficients in $L\cap \CO_{L'}=\CO_L$. (ii): If $m=1$, then $h\circ f\equiv h(X^q)\ (\bmod \fkp)$ proves the claim. Use induction for $m>1$. If $h\circ f = \pi^m g$, then by induction $h=\pi^{m-1}\cdot h'$, thus $h'\circ f = \pi g$ but the $m=1$ case implies $h'\equiv 0\ (\bmod \fkp)$. (iii): If $g(X+_f\alpha)=g(X)$ for all $\alpha\in \MU_1$, then we can write $g(X)-g(0)=g_1(X)\cdot f(X)$ by (i). Now as $f(X+_f\alpha)=f(X)+_{f^\varphi}f(\alpha)=f(X)$, we have $g_1(X+_f\alpha)=g_1(X)$. Repeating this procedure and setting $g_0:=g$ and $g_i(X)-g_i(0)=g_{i+1}(X)\cdot f(X)$, we get $g(X)=\sum_{i=0}^\infty g_i(0)\cdot f(X)^i$, hence $h(X):=\sum_{i=0}^\infty g_i(0)X^i$ gives $g=h\circ f$. Uniqueness follows from (ii), which implies $h\circ f=0\Lra h=0$.
\pfend

\begin{defn}[Coleman \cite{Col}, de Shalit \cite{deSh}]
For $g\in \CO_L[[X]]$, coefficients of the product $\prod_{\alpha\in \sMU_1}g(X+_f\alpha)$ are $\mO_L$-polynomials in the symmetric functions of $\MU_1$, hence they lie in $\CO_L$. Therefore by Lemma \ref{picompimg}(iii), we get a unique $N(g)\in \CO_L[[X]]$ satisfying:
\begin{equation} \label{normdef}
N(g)\circ f(X) = \prod_{\alpha\in \sMU_1}g(X+_f\alpha).
\end{equation}
Clearly $N(g_1g_2)=N(g_1)N(g_2)$. Also, we set $N^0(g):=g$ and
\[N^m(g):= \bigl( N^{m-1}(N(g)^{\varphi^{-1}}) \bigr)^\varphi \ \ (m\geq 1). \]
If we write $N=N_f$ (called the {\em Coleman operator}), this means $N^m=N_{f^{\varphi^{m-1}}} \circ \cdots \circ N_{f^{\varphi}}\circ N_f$.
\end{defn}

\begin{lem} \label{colemanop1}
For $m\geq 1$, we have $N^m(g)\circ f_m(X)=\prod_{\alpha\in \sMU_m}g(X+_f\alpha)$. 
\end{lem}

\pfbegin
The case $m=1$ is the definition. Use induction on $m$. Fix a set $C$ of representatives of $\MU_m/\MU_1$ as $\mO$-modules, and extend $\varphi$ to a $\widetilde{\varphi}\in \Gal(K_x^m/K)$ (Proposition \ref{invartmap}(i)). Then: 
\[\prod_{\alpha\in \sMU_m}g(X+_f\alpha) = \prod_{\beta\in C}\prod_{\alpha\in \sMU_1}g(X+_f\beta+_f\alpha)=\prod_{\beta\in C}N(g)\circ f(X+_f\beta),\]
and $f(X+_f\beta)=f(X)+_{f^\varphi}f(\beta)$, but as $C\ni\beta\mapsto f(\beta)^{{\widetilde{\varphi}}^{-1}}\in \MU_{m-1}$ is a bijection, 
\[\text{RHS }=\prod_{\alpha\in \sMU_{m-1}}N(g)(f(X)+_{f^\varphi} \alpha^{\widetilde{\varphi}})
 = \Bigl(\prod_{\alpha\in \sMU_{m-1}}N(g)^{\varphi^{-1}}(f^{\varphi^{-1}}(X)+_f\alpha)\Bigr)^\varphi\]
equals $\bigl(N^{m-1}\bigl(N(g)^{\varphi^{-1}}\bigr)\circ f_{m-1}(f^{\varphi^{-1}}(X))\bigr)^\varphi = N^m(g)\circ f_m(X)$ by inductive hypothesis.
\pfend

\begin{lem} \label{colemanop2}
\benu
\item $N(g)\equiv g^\varphi\ (\bmod \fkp)$. In particular, $N(\CO_L[[X]]^\times)\subset \CO_L[[X]]^\times$.
\item For $m\geq 1$, if $g\equiv 1\ (\bmod \fkp^m)$, then $N(g)\equiv 1\ (\bmod \fkp^{m+1})$.
\item If $g\in \CO_L[[X]]^\times$ and $m\geq 1$, then $N^m(g)/N^{m-1}(g)^\varphi\equiv 1\ (\bmod \fkp^m)$.
\eenu
\end{lem}

\pfbegin
(i): As $f(X)\equiv X^q\ (\bmod \fkp)$, LHS of (\ref{normdef}) $\equiv N(g)(X^q)\ (\bmod \fkp)$. On the other hand, if we write $L'=K_x^1$, then $\MU_1\subset \fkp_{L'}$, hence $g(X+_f\alpha)\equiv g(X)\ (\bmod \fkp_{L'})$ for all $\alpha\in \MU_1$. Therefore RHS of (\ref{normdef}) $\equiv g(X)^q\equiv g^\varphi(X^q)\ (\bmod \fkp_{L'})$, and we see $N(g)\equiv g^\varphi\ (\bmod \fkp)$. (ii): If we let $g=1+\pi^mh$ and $L'=K_x^1$, then
\begin{align*}
N(g)\circ f &= \prod_{\alpha\in \sMU_1}\bigl(1+\pi^m h(X+_f\alpha)\bigr) \equiv \bigl(1+\pi^m h(X)\bigr)^q\ (\bmod \fkp^m\fkp_{L'})\\
&\equiv 1+q\pi^mh(X)+\cdots \pi^{mq}h(X)^q \equiv 1\ (\bmod \fkp^m\fkp_{L'}),
\end{align*}
hence $(N(g)-1)\circ f \equiv 0\ (\bmod \fkp^m\fkp_{L'})$, and as it belongs to $\CO_L[[X]]$ we have $(N(g)-1)\circ f\equiv 0\ (\bmod \fkp^{m+1})$. Therefore, by Lemma \ref{picompimg}(ii), we get $N(g)-1\equiv 0\ (\bmod \fkp^{m+1})$. (iii): As $N(g)/g^\varphi\equiv 1\ (\bmod \fkp)$ from (i), apply (ii) to this $m-1$ times.
\pfend

\begin{defn}
For a finite separable extension $K'/K$, we denote the image $N_{K'/K}(K^{\prime \times})$ of the norm map $N_{K'/K}:K^{\prime \times}\ra K^\times$ by $N(K'/K)$. For any separable extension $E/K$, define $N(E/K) :=\bigcap_{K'}N(K'/K)$ where $K'$ runs through all the finite extensions in $E$.
\end{defn}

\begin{prop} \label{normgpsmall}
$N(K_x^m/K) = (1+\fkp^m)\times \langle x \rangle$ for all $m\geq 1$. 
\end{prop}

\pfbegin
Write $L'=K_x^m$ and take $\alpha\in \MU_m^\times$. By Proposition \ref{LTextension}(ii) we have $L'^\times=\CO_{L'}^\times\times \langle -\alpha \rangle$ and $N_{L'/K}(-\alpha)=N_{L/K}(\pi^{\varphi^{m-1}})=x$, hence it suffices to show $N_{L'/K}(\CO_{L'}^\times)=1+\fkp^m$. First we show $N_{L'/K}(\CO_{L'}^\times)\subset 1+\fkp^m$. By the following Lemma \ref{compramgen}, any $u\in \CO_{L'}^\times$ can be written as $u=g(\alpha),\ g\in \CO_L[[X]]^\times$. For $i\geq 0$, set $u_i:=N^i(g)(0)$. Then by Lemma \ref{colemanop1} we have $u_i=\prod_{\alpha\in \sMU_i}g(\alpha)$, hence $N_{L'/L}(u)=\prod_{\alpha\in \sMU_m^\times}g(\alpha)=u_m/u_{m-1}$. Lemma \ref{colemanop2}(iii) shows that $u_m/u_{m-1}^\varphi\in 1+\fkp_L^m$. Hence $N_{L'/K}(u)=N_{L/K}(u_m/u_{m-1})=N_{L/K}(u_m/u_{m-1}^\varphi)\in N_{L/K}(1+\fkp^m_L)\subset 1+\fkp^m$. The other inclusion (not used in the sequel) is seen as follows: as $K_x^m$ is the fixed field of $\Art_K((1+\fkp^m)\times \langle x \rangle)$ by Proposition \ref{artchar}, if $x'/x\in 1+\fkp^m$ then $K_x^m=K_{x'}^m$. Therefore $x'\in N(K_{x'}^m/K)=N(K_x^m/K)$ and $1+\fkp^m\subset N(K_x^m/K)$.
\pfend

\begin{lem} \label{compramgen}
If $L'/L$ is totally ramified and $\alpha$ is a uniformizer of $L'$, then $\CO_{L'}=\CO_L[\alpha]$. 
\end{lem}

\pfbegin
If $[L':L]=n$ and $x=\sum_{i=0}^{n-1}a_i\alpha^i\ (a_i\in L)$, then $v_{L'}(x)= \min_i\{v_{L'}(a_i\alpha^i)\}$, as $v_{L'}(a_i\alpha^i)$ are all distinct. Thus (i) $x=0\Rightarrow a_i=0\ (\forall i)$, (ii) $x\in \CO_{L'}\Leftrightarrow a_i\in \CO_L\ (\forall i)$. By (i), the set $\{1,\alpha,\alpha^2\ldots,\alpha^{n-1}\}$ is a basis of $L'$ over $L$. This and (ii) imply $\CO_{L'}\subset \CO_L[\alpha]$. 
\pfend

\bcor \label{ramnorms}
If $E/L$ is totally ramified and $E$ contains $K_x^\ram$, then $N(E/K)=\langle x \rangle$.
\ecor

\bprf
Proposition \ref{normgpsmall} and $\bigcap_{m\geq 1}(1+\fkp^m)=\{1\}$ imply $N(E/K)\subset N(K_x^\ram/K) \subset \langle x \rangle$, and $N(E/K)$ contains an element with valuation $[L:K]$ by the following lemma. \eprf

\blem
Let $P=P_L:=v_L^{-1}(1)$ be the set of all uniformizers of a local field $L$, and $E/L$ a totally ramified extension. Then $N(E/L)^P:=N(E/L)\cap P$ is non-empty.
\elem

\bprf
If $L'/L$ is finite totally ramified, then $N(L'/L)^P\neq \emptyset$ as $N_{L'/L}$ maps $P_{L'}$ into $P$. For a uniformizer $\pi$ of $L$, we have $P=\pi\cdot \mO_L^\times=\invlim{}P/(1+\fkp_L^m)$, where the quotient is taken by the multiplicative action. As $N_{L'/L}(1+\fkp_{L'}^m)\subset 1+\fkp_L^m$ for all $m\geq 1$, the $N_{L'/L}$ is the $\invlim{}$ of $N_m=N_m^{L'}:P_{L'}/(1+\fkp_{L'}^m)\ra P_L/(1+\fkp_L^m)$. We show $N(L'/L)^P=\invlim{}\bigl({\rm Im} N_m\bigr)$ as subsets of $P$. If $\pi=(\pi_m)_m\in \invlim{}\bigl({\rm Im} N_m\bigr)$, then there is $\pi'\in \invlim{}N_m^{-1}(\pi_m)$ as the $\invlim{}$ of non-empty finite sets is non-empty, and $N(\pi')=\pi$. Converse is clear. Now for general $E/L$, every finite $L'/L$ contained in $E$ is totally ramified, and if $L',L''\subset E$ then $L'L''\subset E$ and ${\rm Im} N_m^{L'L''}\subset {\rm Im} N_m^{L'}\cap {\rm Im} N_m^{L''}$. Hence the intersection $\bigcap_{L'}{\rm Im} N_m^{L'}$ in the finite set $P/(1+\fkp_L^m)$, where $L'$ runs through all finite extensions in $E$, is non-empty. Thus $\invlim{}\bigl(\bigcap_{L'}{\rm Im} N_m^{L'}\bigr)\neq \emptyset$, and it is contained in $\invlim{}\bigl({\rm Im} N_m^{L'}\bigr)=N(L'/L)^P$ for all $L'$, hence in $N(E/L)^P$.
\eprf

\subsection{Base change and LCFT for Lubin-Tate extensions}

\bpr \label{ramLTnorm}
For $\sigma\in W(K^{\rm sep}/K)$ with $v(\sigma)>0$, let $E_\sigma\subset K^{\rm sep}$ be its fixed field. Then $N(E_\sigma/K)=\langle \Art^{-1}(\sigma|_{K^\LT}) \rangle$. 
\epr

\pfbegin
Let $x:=\Art^{-1}(\sigma|_{K^\LT})$. By Proposition \ref{artchar}, we have $K^\ram_x\subset E_\sigma$, and $E_\sigma\cap K^\ur$ is the unramified extension of $K$ of degree $v(\sigma)=v(x)$. Hence Corollary \ref{ramnorms} applies.
\pfend

\begin{thm} \label{basechange} {\em (Base change)}
For a finite separable $K'/K$, we have $K^\LT\subset K^{\prime \LT}$ and the following commutes, i.e.\ for all $x'\in K^{\prime \times}$ we have $\Art_{K'}(x')|_{K^\LT}=\Art_K(N_{K'/K}(x'))$. 
\[\xymatrix{
K^{\prime \times} \ar[r]^-{\Art_{K'}}\ar[d]_{N_{K'/K}} & \Gal(K^{\prime \LT}/K') \ar[d]^-\res\\
K^\times \ar[r]^-{\Art_K} & \Gal(K^\LT/K)
}\]
\end{thm}

\pfbegin
Take $x\in \fkp_{K'}\cap K^{\prime \times}$, and extend $\Art_{K'}(x)\in W(K^{\prime \LT}/K')$ to $\sigma\in W(K^{\rm sep}/K')$. By Proposition \ref{ramLTnorm}, we have $\langle N_{K'/K}(x) \rangle = N_{K'/K}(\langle x \rangle) = N_{K'/K}(N(E_\sigma/K')) = N(E_\sigma/K) = \langle \Art_K^{-1}(\sigma|_{K^\LT}) \rangle$. As $v_K(\sigma|_{K^\LT})=f(K'/K)v_{K'}(\sigma)=f(K'/K)v_{K'}(x)=v_K(N_{K'/K}(x))$, we obtain $N_{K'/K}(x) = \Art_K^{-1}(\sigma|_{K^\LT})$. Therefore $\sigma|_{K^\LT}= \Art_K(N_{K'/K}(x))$ depends only on $\sigma|_{K^{\prime LT}}$, which shows $K^\LT\subset K^{\prime \LT}$ and the commutativity, as $\fkp_{K'}\cap K^{\prime \times}$ generates $K^{\prime \times}$.
\eprf

\begin{cor} \label{LCFTminusLKW}
{\em (LCFT minus Local Kronecker-Weber)}
\begin{enumerate}\renewcommand{\labelenumi}{(\roman{enumi})}
\item There is a unique homomorphism $\Art_K:K^\times\ra \Gal(K^\LT/K)$ satisfying:
\begin{enumerate}
\item if $\pi$ is a uniformizer of $K$, then $\Art_K(\pi)|_{K^\ur}=\Frob_K$, and 
\item if $K'/K$ is a Lubin-Tate extension, then $\Art_K(N(K'/K))|_{K'}=\id$. 
\end{enumerate}
Moreover, the $\Art_K$ is an isomorphism onto $W(K^\LT/K)\subset \Gal(K^\LT/K)$.
\item If $K'/K$ is finite separable, then $K^\LT\subset K^{\prime \LT}$, and $\Art_{K'}(x)|_{K^\LT}=\Art_K(N_{K'/K}(x))$ for all $x\in K^{\prime \times}$. The $\Art_K$ induces $K^\times/N(K'/K)\isom \Gal((K'\cap K^\LT)/K)$.
\end{enumerate}
\end{cor}

\pfbegin
(i): The map $\Art_K$ satisfies (a) by definition, and (b) by Theorem \ref{basechange}. Conversely, if $\Art'_K$ satisfies these, then for any uniformizer $\pi$ of $K$, (b) and Proposition \ref{LTextension}(ii) imply $\Art'_K(\pi)|_{K_\pi^\ram}=\id$. This and (a) show $\Art'_K(\pi)=\Art_K(\pi)$ by Proposition \ref{artchar}. As $K^\times$ is generated by the uniformizers, we get $\Art'_K=\Art_K$. The last claim was seen in Definition \ref{artmap}. (ii): The first part is Theorem \ref{basechange}, and $\Art_K$ induces $K^\times/N(K'/K)\cong W(K^\LT/K)/\Image(W(K^{\prime \LT}/K'))$. This is isomorphic to $\Gal((K'\cap K^\LT)/K)$, as $W(K^\LT/K)$ surjects onto $\Gal((K'\cap K^\LT)/K)$ and $W(K^{\prime \LT}/K')$ is the inverse image of $W(K^\LT/K)$ under $\Gal(K^{\prime \LT}/K')\ra \Gal(K^\LT/K)$.
\pfend

Above proof of (i) shows that we only need totally ramified Lubin-Tate extensions for the characterization of $\Art_K$. The classical theorems of LCFT for Lubin-Tate extensions (instead of abelian extensions) follow easily from Corollary \ref{LCFTminusLKW}, for example:

\benu
\item  For any finite $K'/K$, we have $N(K'/K)=N((K'\cap K^\LT)/K))$ and $[K^\times:N(K'/K)]\leq [K':K]$. Equality holds if and only if $K'/K$ is Lubin-Tate.
\item  If $K'/K$ is finite and $K''/K$ is Lubin-Tate, then $N(K'/K)\subset N(K''/K)\iff K''\subset K'$. If both are Lubin-Tate, then $N(K''/K)/N(K'/K)\cong \Gal(K'/K'')$ by $\Art_{K'}$.
\item  If $K',K''/K$ are Lubin-Tate extensions, then: \\
$N(K'K''/K)=N(K'/K)\cap N(K''/K),\ \ \ N((K'\cap K'')/K)=N(K'/K)N(K''/K)$.
\item  ({\em Existence theorem}) For any finite index subgroup $H\subset K^\times$ containing $1+\fkp^m$ for some $m$, there is a unique Lubin-Tate extension $K'/K$ such that $N(K'/K)=H$.
\eenu

\section{The Local Kronecker-Weber theorem}

We finish the proof of Theorem A by proving the local Kronecker-Weber theorem, i.e.\ $K^\LT=K^\ab$. This follows easily from the Hasse-Arf theorem (Gold \cite{Gold} or Iwasawa \cite{Iw}, \S7.4; see also Lubin \cite{Lubin}, Rosen \cite{Rosen}). We first prove the Hasse-Arf theorem following Sen \cite{Sen}.

\subsection{Ramification groups}

Let $K'/K$ be a finite totally ramified Galois extension of local fields, and set $G:=\Gal(K'/K)$. For a uniformizer $\pi$ of $K'$, we have $\CO_{K'}=\mO[\pi]$ by Lemma \ref{compramgen}. We write $v:=v_{K'}$ and $q=|\CO/\fkp|=|\CO_{K'}/\fkp_{K'}|$.

\begin{defn}
Let $i(\sigma):=v(\sigma(\pi)-\pi)$, where we set $i(\id)=\infty$. For $n\geq 0$, define $G_n:=\{\sigma\in G \mid i(\sigma)>n\}=\{ \sigma\in G\mid \sigma(\pi)/\pi\in 1+\fkp_{K'}^n \}$. Then $G=G_0$ as $K'/K$ is totally ramified, and $G_n=\{\id\}$ for sufficiently large $n$. They are normal subgroups of $G$, independent of the choice of $\pi$, because $G_n=\{\sigma\in G\mid v(\sigma(a)-a)>n \text{ for all } a\in \CO_{K'} \}$ is the kernel of the group homomorphism $G\ni \sigma\longmapsto \sigma|_{\CO_{K'}} \bmod \fkp_{K'}^{n+1}\in \Aut(\CO_{K'}/\fkp_{K'}^{n+1})$.
\end{defn}

\begin{prop} \label{filtinj}
For $n\in \Z_{\geq 0}$, we have the following injective group homomorphisms, independent of the choice of $\pi$ (they show that $G$ is supersoluble):
\begin{align*}
\theta_0: G_0/G_1\ni \sigma &\longmapsto \sigma(\pi)/\pi\ \bmod \fkp_{K'}\in (\CO_{K'}/\fkp_{K'})^\times \cong \F_q^\times,\\
\theta_n: G_n/G_{n+1}\ni \sigma &\longmapsto (\sigma(\pi)/\pi)-1\ \bmod \fkp_{K'}^{n+1}\in \fkp_{K'}^n/\fkp_{K'}^{n+1} \cong \F_q\ \ (n\geq 1).
\end{align*}
\end{prop}

\pfbegin
The maps are well-defined and injective by definition of $G_n$. For a different uniformizer $\pi'=u\pi$ with $u\in \CO_{K'}^\times$, we have $\sigma(\pi')/\pi' = (\sigma(\pi)/\pi) \cdot (\sigma(u)/u)$, and if $\sigma\in G_n$ then $\sigma(u) \equiv u\ (\bmod \fkp_{K'}^{n+1})$, hence $\sigma(u)/u\in 1+\fkp_{K'}^{n+1}$, hence the maps $\theta_n$ do not depend on the choice of $\pi$. For $\sigma,\tau\in G_n$, if $u=\tau(\pi)/\pi$, then $\sigma\tau(\pi)/\pi=(\sigma(\pi)/\pi)\cdot (\tau(\pi)/\pi) \cdot (\sigma(u)/u)$, and as $u\in \CO_{K'}^\times$ we have $\sigma(u)/u\in 1+\fkp_{K'}^{n+1}$, therefore $\theta_n$ are group homomorphisms.
\pfend

\begin{cor} \label{tamedivisible}
If $G$ is abelian and $G_n\neq G_{n+1}$, then $e_0:=|G_0/G_1|$ divides $n$.
\end{cor}

\pfbegin
Let $\tau\in G_n$ and $\sigma\in G$. We compute $\theta_n(\sigma \tau \sigma^{-1})$ using $\pi'=\sigma^{-1}(\pi)$. If $\tau(\pi')=\pi'(1+a)$ for $a\in \fkp_{K'}^n$, then $\theta_n(\tau) = a \bmod \fkp_{K'}^{n+1}$ by definition. Then $\sigma\tau\sigma^{-1}(\pi)=\sigma\tau(\pi')=\sigma(\pi'(1+a))=\pi(1+\sigma(a))$, hence $\theta_n(\sigma\tau\sigma^{-1})=\sigma(a) \bmod \fkp_{K'}^{n+1}$. If we write $a=b\pi^n$ for $b\in \CO_{K'}$ and $\sigma(\pi)=u\pi$ for $u\in \CO_{K'}^\times$, then $\sigma(a)=\sigma(b)\sigma(\pi)^n=\sigma(b)u^n\pi^n$, and as $\sigma(b)\equiv b \bmod \fkp_{K'}$, we have $\sigma(a)\equiv bu^n\pi^n=u^na\ (\bmod \fkp_{K'}^{n+1})$. Therefore $\theta_n(\sigma\tau\sigma^{-1}) = u^na \bmod \fkp_{K'}^{n+1}$. If $G$ is abelian, then $\sigma\tau\sigma^{-1}=\tau$, hence $a\equiv u^na \bmod \fkp_{K'}^{n+1}$. If $G_n\neq G_{n+1}$, we can choose $\tau\in G_n$ with $\theta_n(\tau)\neq 0$, i.e.\ $a\in \fkp_{K'}^n\setminus \fkp_{K'}^{n+1}$. Also, choose $\sigma\in G$ which generates $G_0/G_1$, i.e.\ $\theta_0(\sigma)=u\bmod \fkp$ has order $e_0$ in $(\CO_{K'}/\fkp_{K'})^\times$. Then $a\equiv u^na\ (\bmod \fkp_{K'}^{n+1})$ implies $e_0\mid n$.
\pfend

\begin{lem} \label{wildcycram}
For $\sigma\in G_1$, we have $v\bigl(\sum_{i=0}^{p-1}\sigma^i(\alpha)\bigr) > v(\alpha)$ for all $\alpha\in K^{\prime \times}$.
\end{lem}

\begin{proof}
Replacing $\alpha$ by $\alpha x$ for $x\in K^\times$, we can assume $\alpha\in \mO_{K'}$. Let $(\sigma-1)(\alpha):=\sigma(\alpha)-\alpha$. Then $\sigma\in G_1$ implies $v((\sigma-1)^{p-1}(\alpha))>\cdots >v((\sigma-1)(\alpha))>v(\alpha)$. The claim follows by $\sum_{i=0}^{p-1}\sigma^i(\alpha)\equiv (\sigma-1)^{p-1}(\alpha)\ (\bmod p\alpha)$, which follows from $(-1)^i\binom{p-1}{i}\equiv 1\ (\bmod p)$. This is seen from $\sum_{i=0}^{p-1}X^i=(X^p-1)/(X-1)=(X-1)^{p-1}$ in $\F_p[X]$.
\end{proof}

\begin{lem} \label{explem}
Let $\sigma\in G_1$. For each $n\in \Z$, there exists $\alpha\in K^{\prime \times}$ such that $v(\alpha)=n$ and $v(\sigma(\alpha)-\alpha)=n+i(\sigma^n)$. Moreover, any $x\in K^{\prime \times}$ can be written as a sum $x=\sum_{n=v(x)}^\infty x_n$ (see Appendix I) where each $x_n$ satisfies above two properties for $n$ if $x_n\neq 0$.
\end{lem}

\pfbegin
For the first part, if $n\geq 0$, then let $\alpha=\prod_{i=0}^{n-1}\sigma^i(\pi)$ for a uniformizer $\pi$ of $K'$ (set $\alpha=1$ for $n=0$). Then clearly $v(\alpha)=n$, and $\sigma(\alpha)/\alpha=\sigma^n(\pi)/\pi$, thus $v(\sigma(\alpha)-\alpha)=v(\alpha)+v((\sigma(\alpha)/\alpha)-1)=n+i(\sigma^n)$. Also, $\alpha^{-1}$ satisfies the properties for $-n$. For the second part, note that $C:=\{0\}\cup \MU_{q-1}$ is a complete set of representatives for $\CO_{K'}\bmod \fkp_{K'}$, and $\sigma$ acts trivially on $C$ as $C\subset K$. Hence we can write $x=\sum_{n=v(x)}^\infty c_n\alpha_n$ where $c_n\in C$ and $\alpha_n$ is the $\alpha$ we constructed above. Thus $x_n:=c_n\alpha_n$ has the required properties if $c_n\neq 0$.
\pfend

\begin{prop}[Sen \cite{Sen}] \label{Sen}
Let $\sigma\in G_1$, and $|\langle \sigma \rangle|=p^m$ for $m\geq 1$ (by Proposition \ref{filtinj}). Let $H_n:=G_n\cap \langle \sigma \rangle$ for $n\geq 1$ and $i_j:=i(\sigma^{p^j})$ for $j\geq 0$ (and $i_j:=\infty$ for $j\geq m$). Then:
\benu
\item $i_{j-1}<i_j$ if $j\leq m$. Also, $H_n=\langle \sigma^{p^j} \rangle$ if and only if $i_{j-1}\leq n < i_j$.
\item $i(\sigma^a)=i_{v_p(a)}$ for $a\geq 1$, where $v_p:=v_{\Q_p}$.
\item $i_{j-1}\equiv i_j\ (\bmod p^j)$, where $\infty$ is understood to be congruent to any integer.
\eenu
\end{prop}

\begin{proof}
(i): Lemma \ref{wildcycram} for $\alpha=\sigma^{p^{j-1}}(\pi)-\pi$ shows $i_{j-1}<i_j$. We have $\langle \sigma^{p^j} \rangle\subset H_n$ if and only if $\sigma^{p^j}\in H_n$, i.e.\ $i_j>n$. As all subgroups of $\langle \sigma \rangle$ are of the form $\langle \sigma^{p^j} \rangle$, we have $\langle \sigma^{p^j} \rangle \supset H_n \Leftrightarrow \langle \sigma^{p^{j-1}} \rangle \not{\!\!\subset}\, H_n \Leftrightarrow i_{j-1}\leq n$. (ii): This is $\infty=\infty$ if $p^m\mid a$. If $j:=v_p(a)<m$, then $H_{i_j-1}=\langle \sigma^{p^j} \rangle$ and $H_{i_j}=\langle \sigma^{p^{j+1}} \rangle$ by (i), therefore $\sigma^a\in H_{i_j-1}\setminus H_{i_j}$, i.e.\ $i(\sigma^a)=i_j$. (iii): We can assume $i_j<\infty$, and use induction on $j$. The assertion is empty when $j=0$. Let $j=1$, and assume the {\em Inductive Hypothesis} (the assertion of (iii) for $j-1$). We first prove the {\em Claim: the $i_{j-1}$ and $n+i(\sigma^n)$ for $n\in \Z,\ v_p(n)<j$ are all distinct from each other}. As $v_p(n)\leq j-1$, the {\em Inductive Hypothesis} shows $i(\sigma^n)=i_{v_p(n)}\equiv i_{j-1}\ (\bmod p^{v_p(n)+1})$, i.e.\ $v_p(i_{j-1}-i(\sigma^n))>v_p(n)$, hence $i_{j-1}\neq n+i(\sigma^n)$. Now assume $n+i(\sigma^n)=n'+i(\sigma^{n'})$. If $v_p(n)\neq v_p(n')$, then $v_p(n-n')=\min\{ v_p(n), v_p(n') \}$, but the {\em Inductive Hypothesis} shows $v_p(i(\sigma^n)-i(\sigma^{n'})) > \min \{v_p(n), v_p(n') \}$, which is impossible. Hence $v_p(n)=v_p(n')$, therefore $i(\sigma^n)=i(\sigma^{n'})$ and $n=n'$. Thus the {\em Claim} is proven. Now applying the {\em Inductive Hypothesis} to $\sigma^p\in G_1$, we have $i_{j-1}\equiv i_j\ (\bmod p^{j-1})$. Let $s:=i_{j-1}-i_j$ and assume $v_p(s)=j-1$, to see it leads to contradiction. The first part of Lemma \ref{explem} for $\sigma^p$ shows that there is $x\in K^{\prime \times}$ with $v(x)=s$ and $v(\sigma^p(x)-x)=s+i((\sigma^p)^s)=s+i_j=i_{j-1}$. Letting $y:=\sum_{i=0}^{p-1}\sigma^i(x)$, we have $v(y)>v(x)=s$ by Lemma \ref{wildcycram} and $v(\sigma(y)-y)=v(\sigma^p(x)-x)=i_{j-1}$. Now expand $y=\sum_{n=v(y)}^\infty y_n$ as in Lemma \ref{explem}: $v(\sigma(y_n)-y_n)=n+i(\sigma^n)$ if $y_n\neq 0$. Let $z:=\sigma(y)-y$. Then $v(z)=i_{j-1}$ and $z=\sum_{n=v(y)}^\infty z_n$, where $z_n:=\sigma(y_n)-y_n$, hence $v(z_n)=n+i(\sigma^n)$ whenever $z_n\neq 0$. The {\em Claim} shows $v\bigl(z-\sum_{v_p(n)<j}z_n\bigr)\leq i_{j-1}$. If $v_p(n)\geq j$ and $z_n\neq 0$, then $v(z_n)=n+i(\sigma^n)\geq n+i_j\geq v(y)+i_j>i_{j-1}$, hence $v\bigl(\sum_{v_p(n)\geq j}z_n\bigr)>i_{j-1}$, a contradiction.
\pfend

\begin{cor} \label{cyclicHA}
Assume $G\cong \Z/p^m\Z$. Then there exist $n_0,n_1,\ldots,n_{m-1}\in \Z_{\geq 1}$ such that, for $1\leq j\leq m-1$, we have $|G_n|=p^{m-j}$ if and only if $\sum_{i=0}^{j-1}n_ip^i < n\leq \sum_{i=0}^{j}n_ip^i$.
\end{cor}

\subsection{The Hasse-Arf theorem}

Let $G=\Gal(K'/K)$ with $K'/K$ totally ramified as before, and let $G\rhd H$ with $G/H=\Gal(K''/K)$. For $\sigma\in G$, let $\ol{\sigma}=\sigma H\in G/H$ be its image.

\begin{lem} \label{iofquot}
For all $\sigma\in G$, we have $i(\ol{\sigma})=\frac{1}{|H|}\sum_{\tau\in H}i(\sigma \tau)$.
\end{lem}

\pfbegin
For $\ol{\sigma}=\id$, we understand the equality as $\infty=\infty$. Let $\ol{\sigma}\neq \id$, and take uniformizers $\pi'$ and $\pi''$ of $K'$ and $K''$ respectively, so that $\CO_{K'}=\CO[\pi']$ and $\CO_{K''}=\CO[\pi'']$ by Lemma \ref{compramgen}. As $i(\ol{\sigma})=v_{K''}(\ol{\sigma}(\pi'')-\pi'')=\frac{1}{|H|}\cdot v_{K'}(\ol{\sigma}(\pi'')-\pi'')$, if we let $a=\ol{\sigma}(\pi'')-\pi''$ and $b=\prod_{\tau\in H}(\sigma\tau(\pi')-\pi')$, it suffices to show $v_{K'}(a)=v_{K'}(b)$. Let the minimal polynomial of $\pi'$ over $\CO_{K''}$ be $f=\prod_{\tau\in H}(X-\tau(\pi'))\in \CO_{K''}[X]$. Applying $\sigma$, we get $f^{\ol{\sigma}}=\prod_{\tau\in H}(X-\sigma\tau(\pi'))$, where $f^{\ol{\sigma}}\in \CO_{K''}[X]$ is obtained by applying $\ol{\sigma}$ to the coefficients of $f$. Hence $f^{\ol{\sigma}}(\pi')=\prod_{\tau\in H}(\pi'-\sigma\tau(\pi'))=\pm b$. First we prove $a\mid b$. As $\CO_{K''}=\CO[\pi'']$, we have $a\mid \ol{\sigma}(x)-x$ for any $x\in \CO_{K''}$, hence $a\mid f^{\ol{\sigma}}-f$, therefore $a\mid f^{\ol{\sigma}}(\pi')-f(\pi')=\pm b$. Now we prove $b\mid a$. Write $\pi''=g(\pi')$ for $g\in \CO[X]$. The polynomial $g(X)-\pi''\in \CO_{K''}[X]$ has $\pi'$ as a root, hence divisible by $f$ in $\CO_{K''}[X]$. Applying $\ol{\sigma}$, we have $f^{\ol{\sigma}}\mid g(X)-\ol{\sigma}(\pi'')$ in $\CO_{K''}[X]$, hence $g(\pi')-\ol{\sigma}(\pi'')=-a$ is divisible by $f^{\ol{\sigma}}(\pi')=\pm b$.
\pfend

\begin{prop}[Herbrand] \label{Herbrand}
Define $\phi_H(n):=-1+\frac{1}{|H|}\sum_{\tau\in H}\min\{i(\tau), n+1\}$ for $n\in \R_{\geq 0}$.  Also, for $n\in \R_{\geq 0}$, define $G_n:=\{\sigma\in G\mid i(\sigma)\geq n+1\}$, i.e.\ $G_n=G_i$ if $i\in \Z_{\geq 0}$ and $n\in (i-1,i]$. Then $G_nH/H=(G/H)_{\phi_H(n)}$ for all $n\in \R_{\geq 0}$.
\end{prop}

\pfbegin
For $\ol{\sigma}\in G/H$, replace $\sigma$ by the element in $\sigma H$ which has the maximal value of $i$, and let $i(\sigma)=m$. Let $\tau\in H$. If $i(\tau)\geq m$, then $i(\sigma\tau)\geq m$, hence $i(\sigma\tau)=m$. If $i(\tau)<m$, then $i(\tau)\geq \min\{ i(\sigma\tau), i(\sigma^{-1}) \}$, hence $i(\sigma\tau)=i(\tau)$. Therefore $i(\sigma\tau)=\min\{i(\tau), m\}$. Now the Lemma \ref{iofquot} gives $i(\ol{\sigma})=\phi_H(m-1)+1$. Therefore, as $\phi_H$ is increasing, for $n\in \R_{\geq 0}$ we have $\ol{\sigma}\in G_nH/H\iff m\geq n+1 \iff i(\ol{\sigma})\geq \phi_H(n)+1 \iff \ol{\sigma}\in (G/H)_{\phi_H(n)}$.
\pfend

\begin{lem} \label{phitrans}
Let $\phi_G(n):=-1+\frac{1}{|G|}\sum_{\tau\in G}\min\{i(\tau), n+1\}$ for $n\in \R_{\geq 0}$. Then:
\benu
\item $\phi_G(0)=0,\ \phi_G(n)=\frac{1}{|G|}\sum_{i=1}^n|G_i|$ for $n\in \Z_{\geq 1}$. 
\item $\phi_G=\phi_{G/H}\circ \phi_H$ on $\R_{\geq 0}$.
\eenu
\end{lem}

\pfbegin
(i): $\displaystyle \sum_{\tau\in G}\min\{i(\tau), n+1\}=\sum_{i=0}^{n-1}\Bigl(\sum_{\tau\in G_i\setminus G_{i+1}}(i+1) \Bigr) + \sum_{\tau\in G_n}(n+1) = \sum_{i=0}^n|G_i|$.

(ii): As $\phi(0)=0$ and $\phi$ is continuous and piecewise linear, we only need to compare the derivatives of both sides at $n\in (i-1,i)$ for $i\in \Z_{>0}$. For LHS it is $|G_n|/|G|$, and for RHS it is $(|(G/H)_{\phi_H(n)}|/|G/H|)\cdot (|H_n|/|H|)=|G_nH/H||H_n|/|G|=|G_n|/|G|$ by Proposition \ref{Herbrand} and $G_n/H_n=G_n/(H\cap G_n)\cong G_nH/H$.
\pfend

\begin{thm}[Hasse-Arf] \label{HasseArf}
If $G$ is abelian, $n\in \Z_{\geq 0}$ and $G_n\neq G_{n+1}$, then $\phi_G(n)\in \Z_{\geq 0}$.
\end{thm}

\pfbegin
First assume $G=G_1$. Then $G\cong \bigoplus_{i=1}^j\Z/p^{m_i}\Z$ by Proposition \ref{filtinj}, and we proceed by induction on $j$. When $j=1$, i.e.\ $G\cong \Z/p^m\Z$, if $G_n\neq G_{n+1}$ then $n=\sum_{i=0}^{j}n_ip^i$ for some $0\leq j\leq m-1$ by Corollary \ref{cyclicHA}, in which case $\phi_G(n)=\frac{1}{p^m}(n_0\cdot p^m + n_1p\cdot p^{m-1} + \cdots + n_jp^j\cdot p^{m-j})\in \Z_{\geq 0}$ by Lemma \ref{phitrans}(i). For $j>1$, if $G_n\neq G_{n+1}$ we can find $H$ with $G/H\cong \Z/p^{m_i}\Z$, and $G_nH/H\neq G_{n+1}H/H$. We have $\phi_H(n)\in \Z_{\geq 0}$ by inductive hypothesis, and $(G/H)_{\phi_H(n)}\neq (G/H)_{\phi_H(n+1)}=(G/H)_{\phi_H(n)+1}$ by Proposition \ref{Herbrand}. As $G/H$ is cyclic, we see $\phi_{G/H}(\phi_H(n))\in \Z_{\geq 0}$, which is $\phi_G(n)$ by Lemma \ref{phitrans}(ii). Now when $G\neq G_1$, set $H=G_1$ and $|G/H|=e_0$. As $\phi_{G/H}(n)=n/e_0$ for $n\in \R_{\geq 0}$ by definition, by Lemma \ref{phitrans}(ii) it suffices to show $e_0\mid \phi_H(n)$ when $n\in \Z_{\geq 0}$ and $G_n\neq G_{n+1}$ (we know $\phi_H(n)\in \Z_{\geq 0}$). If $n=0$ then $\phi_H(0)=0$. Let $n>0$. For any $i\in \Z_{\geq 1}$ (where $H_i=G_i$) with $H_i\neq H_{i+1}$, we have $e_0\mid i$ by Corollary \ref{tamedivisible}, hence $e_0\mid \sum_{i=1}^n|H_i|$. As $e_0$ and $|H|$ are coprime, we have $e_0\mid \phi_H(n)$ by Lemma \ref{phitrans}(i).
\pfend

\begin{defn} \label{upperno}
For $m\in \R_{\geq 0}$, set $G^m:=G_{\phi_G^{-1}(m)}$ (the {\em upper numbering}). 
\end{defn}

\begin{cor} \label{ramdegree}
\benu
\item If $G\rhd H$, then $G^mH/H=(G/H)^m$ for all $m\in \R_{\geq 0}$.
\item Let $K'/K$ and $K''/K$ be two Galois extensions with $K'K''/K$ totally ramified. If $\Gal(K'/K)^m=\Gal(K''/K)^m=\{\id\}$ for $m\in \R_{\geq 0}$, then $\Gal(K'K''/K)^m=\{\id\}$.
\item Let $G$ be abelian. Then $|G/G^m|$ divides $(q-1)q^{m-1}$ for $m\in \Z_{\geq 0}$.
\eenu
\end{cor}

\pfbegin
(i): By Proposition \ref{Herbrand} and Lemma \ref{phitrans}(ii), we compute $G^mH/H=G_{\phi^{-1}_G(m)}H/H=(G/H)_{\phi_H(\phi^{-1}_G(m))}=(G/H)_{\phi^{-1}_{G/H}(m)}=(G/H)^m$. (ii): If $G=\Gal(K'K''/K)$ and $G/H=\Gal(K''/K)$, then $G^mH/H=(G/H)^m=\{\id\}$ shows $G^m\subset H=\Gal(K'K''/K'')$. Similarly $G^m\subset \Gal(K'K''/K')$, hence $G^m=\{\id\}$. (iii): If $n-1< \phi_G^{-1}(m)\leq n$ for $n\in \Z_{\geq 0}$, then $G^m=G_n$. Consider $G_i$ for integers $1 \leq i\leq n$. Then, by Theorem \ref{HasseArf}, $G_{i-1}\neq G_i$ can only happen when $\phi_G(i-1)\in \Z$, and as $0\leq \phi_G(i-1)\leq \phi_G(n-1)<m$, at most $m-1$ times for $i>1$. By Proposition \ref{filtinj}, $|G_{i-1}/G_i|$ divides $q-1$ when $i=1$ and $q$ when $i>1$.
\pfend

\subsection{The Local Kronecker-Weber theorem}

\begin{prop} \label{LTramnumber}
Let $x\in K^\times$ with $v(x)=n>0$. Let $L=K_n$ and $K_x^m$ as in Definition \ref{finLTkxm}. Then we have $\Gal(K_x^m/L)^m=\{\id\}$ for all $m\geq 1$ (see Definition \ref{upperno}).
\end{prop}

\pfbegin
Let $K_x^m=L_f^m$ and $\alpha\in \MU_{f,m}^\times$. For $\sigma\in \Gal(K_x^m/L)\setminus \{\id\}$, we have $i(\sigma)=v(\sigma(\alpha)-\alpha)$ by Proposition \ref{LTextension}(ii), where $v=v_{K_x^m}$. If $\rho_{f,m}(\sigma)= u\bmod \fkp^m\in (\CO/\fkp^m)^\times$ (see Proposition \ref{LTextension}(iii)), then $\sigma(\alpha)=[u]_f(\alpha)$. For $\sigma\neq \id$, set $\beta:=[u-1]_f(\alpha)$. If $v_K(u-1)=i$ for $0\leq i<m$, then $\beta\in \MU_{f,m-i}^\times$ by Lemma \ref{mupmtors}(ii). Hence $\beta$ is a uniformizer of $K_x^{m-i}$ by Proposition \ref{LTextension}(ii), which shows $v(\beta)=q^i$. Now $\sigma(\alpha)=[u]_f(\alpha)=\alpha +_f \beta \equiv \alpha + \beta\ (\bmod \alpha\beta)$, hence $i(\sigma)=v(\sigma(\alpha)-\alpha)=v(\beta)=q^i$. Thus for $G=\Gal(K_x^m/L)$ and $1\leq i\leq m$, we have $|G_n|=|\rho_{f,m}^{-1}(1+\fkp^i)|=q^{m-i}$ for $q^{i-1}-1 < n \leq q^i-1$. Thus $\phi_G(q^m-1)=\frac{1}{|G|}\sum_{i=1}^{q^m-1}|G_i|=\frac{1}{(q-1)q^{m-1}}\bigl(\sum_{i=1}^m (q^i-q^{i-1})\cdot q^{m-i}\bigr)=m$ and $G^m=G_{q^m-1}=\{\id\}$.
\pfend

\begin{thm} {\em (Local Kronecker-Weber theorem)} 
Every finite abelian extension of a local field $K$ is a Lubin-Tate extension, i.e.\ $K^\LT=K^\ab$.
\end{thm}

\pfbegin
Take a $\sigma\in W(K^\LT/K)$ with $v(\sigma)=n>0$, and let $L=K_n$. Extend $\sigma$ arbitrarily to $\sigma\in W(K^\ab/K)$, and let $E_\sigma\subset K^\ab$ be its fixed field. Then $E_\sigma\cap K^\ur=L$ and $E_\sigma/L$ is totally ramified Galois. Now $\Gal(K^\ab/E_\sigma)\cong \widehat{\Z}$ with $\sigma\mapsto 1$ by the definition of $E_\sigma$. On the other hand, $\Gal(K^\ur E_\sigma/E_\sigma)\cong \Gal(K^\ur/L)\cong \hat{\Z}$ by $\sigma\mapsto 1$, as $\sigma|_{K^\ur}=\Frob_L$. Therefore $\Gal(K^\ab/E_\sigma)\cong \Gal(K^\ur E_\sigma/E_\sigma)$, i.e.\ $K^\ab=K^\ur E_\sigma$. Now set $x:=\Art_K^{-1}(\sigma)$. Then $K_x^\ram\subset E_\sigma$ by Proposition \ref{artchar}. As $K^\LT=K^\ur K_x^\ram $, it suffices to show $E_\sigma\subset K_x^\ram$. Let $K'/L$ be any finite Galois extension contained in $E_\sigma$. It is totally ramified, and $\Gal(K'/L)^m=\{\id\}$ for a large $m$. Then we have $\Gal(K'K_x^m/L)^m=\{\id\}$ by Proposition \ref{LTramnumber} and Corollary \ref{ramdegree}(ii), hence $[K'K_x^m:L]\mid (q-1)q^{m-1} = [K_x^m:L]$ by Corollary \ref{ramdegree}(iii), thus $K'\subset K_x^m$.
\pfend



\renewcommand\thesection{\Roman{section}}

\section*{Appendix I: Basic facts on DVR}
\setcounter{section}{1}
\setcounter{thm}{0}

Here we gather some facts on DVR that are used in this article. The proofs omitted here can be found in Atiyah-Macdonald \cite{AM} and the Chapters I, II of Serre \cite{Serre}. A ring $A$ is called a {\em discrete valuation ring (DVR)} if it is a local ring (i.e.\ has a unique maximal ideal), a PID and not a field. Let $A$ be a DVR with the maximal ideal $P$, and let $K$ be its fraction field. A generator of $P$ is called a {\em uniformizer} of $A$. Each uniformizer $\pi$ gives a following isomorphism of abelian groups:
\[A^\times\times \Z\ni (u,b) \stackrel{\cong}{\longmapsto} u\cdot \pi^b\in K^\times.\]
The second projection ({\em valuation}) $v_K:K^\times\ra \Z$ does not depend on $\pi$, and setting $v_K(0):=\infty$, we have $A=\{x\in K\mid v_K(x)\geq 0\}$ and $P=\{x\in K\mid v_K(x)> 0\}$.

The {\em completion} of $A$ is defined as $\widehat{A}:=\invlim{m}A/P^m$, which is also a DVR with the maximal ideal $\widehat{P}:=P\widehat{A}$. If $K=\Frac(A)$, then $\widehat{K}:=K\otimes_A\widehat{A}$ is the fraction field of $\widehat{A}$, which is called the {\em completion} of $K$. The canonical map $A\ra \widehat{A}$ is always injective (hence $K\subset \widehat{K}$), and if it is an isomorphism we call $A$ a {\em complete discrete valuation ring (CDVR)}. For example, the {\em ring of $p$-adic integers} $\Z_p:=\invlim{}\Z/(p^m)$ is a CDVR with $(p)$ as its maximal ideal, and its fraction field $\Q_p$ is the {\em $p$-adic field}. A completion of a DVR is a CDVR, and $A/P^m\isom \widehat{A}/\widehat{P}^m$. If $A$ is a DVR, choosing a complete set of representatives $C$ for $A \bmod P$ and elements $x_n\in A$ with $v(x_n)=n$ for all $n\geq 0$, we can write any element of $\widehat{A}=\invlim{}A/P^m$ uniquely as $\bigl(\sum_{n=0}^{m-1} c_nx_n \bmod P^m\bigr)_m$ with $c_n\in C$. (Incidentally, this shows that if $|C|<\infty$ then $|A/P^m|=|C|^m$.) We write this element as $\sum_{n=0}^\infty c_nx_n$ (when $x_n=\pi^n$ for a uniformizer $\pi$, this is called a {\em $\pi$-adic expansion}). Choosing $x_n\in K$ with $v(x_n)=n$ for all $n\in \Z$, any $x\in \widehat{K}$ can be written as $y+\sum_{v(x)\leq n<0}c_nx_n$ for some $y\in \widehat{A}$, hence as $x=\sum_{n=v(x)}^\infty c_nx_n$. 

If $A$ is a CDVR, then we can substitute $x_1,\ldots,x_n\in P$ into any power series $F\in A[[X_1,\ldots,X_n]]$ with coefficients in $A$ to get $F(x_1,\ldots,x_n)\in A$. This is defined using $A[[X_1,\ldots,X_n]]\cong \invlim{m}\bigl(A[X_1,\ldots,X_n]/(\deg m)\bigr)$ and $A\cong \invlim{m}A/P^m$, by taking the limit of:
\[A[X_1,\ldots,X_n]/(\deg m)\ni F\bmod \deg m\longmapsto F(x_1,\ldots,x_n)\bmod P^m\in A/P^m.\]

Let $A$ be a DVR, $K$ its fraction field, $L$ a separable extension of $K$ of degree $n$, and $B$ the integral closure of $A$ in $L$, so that $L\cong B\otimes_AK$ and $L=\Frac(B)$. Then $B$ is a finitely generated $A$-module, and as $A$ is a PID, it is a free $A$-module of rank $n=[L:K]$. Also, $B$ is a {\em Dedekind domain}, i.e.\ 1-dimensional integrally closed noetherian domain. If $PB=\prod_{i=1}^gQ_i^{e_i}$ is the prime ideal decomposition of the ideal $PB$ of $B$ generated by the elements of $P$, then $Q_1,\ldots, Q_g$ are all the maximal ideals of $B$. Let $\widehat{B}_i:=\invlim{m}B/Q_i^m$ for $1\leq i\leq g$. As $B$ is a finite free $A$-module, the functor $B\otimes_A$ and inverse limits commute, hence the following canonical maps are isomorphisms: 
\[B\otimes_A\widehat{A}\cong B\otimes \Bigl(\invlim{m}A/P^m\Bigr)\cong \invlim{m}B/(PB)^m\cong \invlim{m}\prod_{i=1}^gB/Q_i^{e_im}\cong \prod_{i=1}^g\widehat{B}_i.\]

\begin{prop} \label{dvrlemma}
\benu
\item If $A$ is a CDVR, then so is $B$.
\item If $B$ is also a DVR, then the completion $\widehat{L}$ of $L$ is isomorphic to $L\otimes_K\widehat{K}$ (i.e.\ it is the composite field $L\widehat{K}$), and $L\cap \widehat{K}=K$ in $\widehat{L}$. 
\eenu
\end{prop}

\pfbegin
(i): $B\cong B\otimes_A\widehat{A}$ and $B$ is a domain, hence $g=1$ and $B\cong \widehat{B}$. (ii): $B\otimes_A\widehat{A}\cong \widehat{B}$ gives $L\otimes_K\widehat{K}\cong L\otimes_K(K\otimes_A\widehat{A})\cong L\otimes_B (B\otimes_A\widehat{A})\cong L\otimes_B\widehat{B}\cong \widehat{L}$. Now let $K':=L\cap \hat{K}$ and $[K':K]=m$. As $K'/K$ is separable, let $K'\cong K[X]/(f)$ with $\deg f=m$. Assume $m>1$. As $f$ has a root in $K'\subset L$, we have $L\otimes_KK'\cong L[X]/(f)\cong L\times L'$ with an $L$-algebra $L'$; but then $\hat{L}\cong L\otimes_K\hat{K}\cong (L\otimes_K K')\otimes_{K'}\hat{K}\cong (L\times L')\otimes_{K'}\hat{K}\cong (L\otimes_{K'}\hat{K})\times (L'\otimes_{K'}\hat{K})$, a contradiction because $\hat{L}$ is a field.
\pfend

Assume $g=1$ and $Q=Q_1$ in the following. As $Q\cap A=P$, the field $k_Q:=B/Q$ is an extension of $k_P:=A/P$, and as $B$ is a finitely generated $A$-module, $k_Q/k_P$ is finite. The {\em ramification index} $e$ and {\em residue degree} $f$ are defined by $PB=Q^e$ and $f=[k_Q:k_P]$. As vector spaces over $k_P$, we have $B/PB\cong (k_Q)^e$ (use $Q$-adic expansion), and the dimension of RHS is $ef$, and the dimension of LHS is the rank of $B$ as an $A$-module, which is $n$. Therefore $n=ef$. Assume moreover that $L/K$ is Galois and $k_P$ is perfect. We say $L/K$ is {\em unramified} if $e=1$ and {\em totally ramified} when $f=1$. An element of $\Gal(L/K)$ induces an automorphism of $B$ which maps $Q$ onto itself, hence we have a group homomorphism: 
\[ \Gal(L/K)\ni \sigma\longmapsto \sigma|_B\bmod Q\in \Aut(k_Q/k_P). \]
We can show that $k_Q/k_P$ is Galois and the homomorphism is surjective. As $|\Gal(k_Q/k_P)|=f$, the order of the kernel is $e$. The following gives an unramified example: 

\begin{prop} \label{cyclounram}
Let $L=K(\MU_n)$ (and $g=1$). If $\chara k_P\notmid n$, then $L/K$ is unramified. 
\end{prop}

\pfbegin
We show that the above homomorphism is injective. As any element of $\Gal(K(\MU_n)/K)$ is determined by the image of a generator $\zeta\in B^\times$ of $\MU_n$, it suffices to show that if $\zeta^i\equiv \zeta^j\ (\bmod Q)$ then $\zeta^i=\zeta^j$. As $\zeta^i-\zeta^j\in Q$ implies $\zeta^{i-j}-1\in Q$, we only need to show $\zeta^i-1\notin Q$ for $1\leq i\leq n-1$. Substituting $X=1$ to the identity $\prod_{i=1}^{n-1}(X-\zeta^i)=(X^n-1)/(X-1)=X^{n-1}+X^{n-2}+\cdots +X+1$, we get $\prod_{i=1}^{n-1}(1-\zeta^i)=n$, and as $n\notin Q$ we have $\prod_{i=1}^{n-1}(1-\overline{\zeta^i})\neq 0$ in the field $k_Q$, hence $\zeta^i-1\notin Q$. 
\pfend

\begin{proof}[Proof of Lemma \ref{Hensel}]
$(\Leftarrow)$ follows from $\zeta^i\equiv \zeta^j\ (\bmod \fkp)\Lra \zeta^i=\zeta^j$, which we showed in the proof of Proposition \ref{cyclounram}. We show $(\Rightarrow)$. As there is a generator of $\MU_n$ in $k=\CO/\fkp$, take its representative $\zeta_1\in \CO$. As $\CO=\invlim{} \CO/\fkp^m$, it is enough to construct $\zeta_m\in \CO$ for each $m\geq 1$ such that $\zeta_m^n\equiv 1\ (\bmod \fkp^m)$ and $\zeta_{m+1}\equiv \zeta_m\ (\bmod \fkp^m)$. If we have $\zeta_m$, let $\zeta_m^n\equiv 1+\alpha\pi^m\ (\bmod \fkp^{m+1})$. Setting $\zeta_{m+1}=\zeta_m+\beta\pi^m$, we need $\zeta_{m+1}^n\equiv \zeta_m^n+n\zeta_m^{n-1}\beta\pi^m\equiv 1+(\alpha+n\zeta_m^{n-1}\beta)\pi^m\ (\bmod \fkp^{m+1})$ to be $\equiv 1\bmod \fkp^{m+1}$, hence $\beta=-\alpha/n\zeta_m^{n-1}$ will do.
\pfend

\section*{Appendix II: Separability of $f_m$} \label{separability}
\setcounter{section}{2}
\setcounter{thm}{0}

Here we prove the separability of $f_m$ of Definition \ref{fmdef} directly. It is used in the proof of Lemma \ref{mupmtors}(i) only when $\chara K=p$. On the other hand, it follows from Proposition \ref{LTextension}(i) when $\chara K=0$.

\bpr \label{fmsep}
For $\forall m\geq 0$, $f_m\in \mO_L[X]$ is separable.
\epr

\bprf
Lemma \ref{fmseplem2} will show that $f_m(\alpha)=0\Lra f'_m(\alpha)\neq 0$ for all $\alpha\in \ol{L}$.
\eprf

\blem \label{fmseplem1}
Let $\mO'$ be an $\mO_L$-algebra, and $f\in \mO_L[X]$ as above.
\benu
\item Let $\mO'$ be a domain and $\alpha\in \mO'$. If $\alpha\notin \mO^{\prime \times}$, then $f'(\alpha)\neq 0$.
\item Let $\mO'$ be a domain and integral over $\mO_L$, and $f(\alpha)=\beta$ for $\alpha,\beta\in \mO'$. If $\alpha\in \mO^{\prime \times}$, then (a) $\beta\neq 0$, and (b) if $\beta\mid \pi$ in $\mO'$, then $\beta\in \mO^{\prime \times}$.
\eenu
\elem

\bprf
(i): As $\pi\mid q$ in $\mO$, we have $f'(X)=\pi(1+X g(X))$ with $g\in \mO_L[X]$, hence if $\alpha\notin \mO^{\prime \times}$, then $1+\alpha g(\alpha)\neq 0$ and $f'(\alpha)\neq 0$. (ii): As $\beta=f(\alpha)=\alpha^n+\pi g(\alpha)$ with $g\in \mO_L[X]$, we have $\beta-\pi g(\alpha)\in \mO^{\prime \times}$ if $\alpha\in \mO^{\prime \times}$. As $\pi\notin \mO^{\prime \times}$ because $\mO_L$ is integrally closed, we have $\beta\neq 0$. If $\pi=\beta\beta'$, then $\beta(1-\beta'g(\alpha))\in \mO^{\prime \times}$, hence $\beta \in \mO^{\prime \times}$.
\eprf

\blem \label{fmseplem2}
Let $\alpha\in \ol{L}$, and let $\mO_L[\alpha]$ be the $\mO_L$-subalgebra of $\ol{L}$ generated by $\alpha$.
\benu
\item If $f_i(\alpha)\notin \mO_L[\alpha]^\times$ for all $0\leq i\leq m-1$, then $f'_m(\alpha)\neq 0$.
\item If $f_m(\alpha)=0$, then $f_i(\alpha)\notin \mO_L[\alpha]^\times$ for all $0\leq i\leq m-1$.
\eenu
\elem

\bprf
(i): The claim is empty when $m=0$ as $f_0'(X)=1$. We prove by induction on $m$: assume it is true for $m-1$. As $f_{m-1}(\alpha)\notin \mO[\alpha]^\times$, by Lemma \ref{fmseplem1}(i), we have $(f^{\varphi^{m-1}})'(f_{m-1}(\alpha))\neq 0$. By the induction hypothesis, we have $f'_{m-1}(\alpha)\neq 0$. Hence $f'_m(\alpha)\neq 0$. (ii): If $f_i(\alpha)=0$, then $f_j(\alpha)=0$ for $\forall j\geq i$, so we can assume $f_i(\alpha)\neq 0$ for $0\leq i\leq m-1$. Then we have $\alpha\mid f(\alpha)\mid \cdots \mid f_{m-1}(\alpha)\mid \pi^{\varphi^{m-1}}$ in $\mO[\alpha]$, which is finite, hence integral, over $\mO$, as $\alpha$ is a root of a monic $f_m\in \mO_L[X]$. Now assume $f_i(\alpha)\in \mO[\alpha]^\times$ for some $i$. If $i\neq m-1$, then $f_{i+1}(\alpha)\mid \pi$, hence $f_{i+1}(\alpha)\in \mO[\alpha]^\times$ by Lemma \ref{fmseplem1}(ii). Therefore $f_{m-1}(\alpha)\in \mO[\alpha]^\times$, but then $f_m(\alpha)\neq 0$ by Lemma \ref{fmseplem1}(ii), a contradiction.
\eprf

\section*{Remarks on the literature} \label{literature}

The ``relative" Lubin-Tate groups treated in \S3, \S4 and \S5 are due to de Shalit \cite{deSh}, although proofs are omitted there. The exposition is based on Iwasawa \cite{Iw}, with two notable differences. Firstly, in Iwasawa \cite{Iw} the norm operator $N$ is treated only for the ``classical" Lubin-Tate groups, which proves the base change theorem for totally ramified extensions (and the part (i) of Theorem A), and then appeals to the local Kronecker-Weber theorem to prove the base change in the unramified case. Here we provided a uniform proof by using the norm operator in the general setting. Secondly, we separated the ``geometric" (\S3, \S4) and ``arithmetic" (\S5) parts of the theory by defining the Artin map through an arbitrary Lubin-Tate group over $\hat{\mO}$, in the spirit of Carayol \cite{Car}. In \S6 we combined Sen \cite{Sen} with the standard material from Serre \cite{Serre}, Chapter IV. Throughout this article we avoided the use of topological rings/fields, and instead used the language of commutative algebra, which might be a somewhat new way of exposition. Needless to say, there are many other important approaches to local class field theory, see e.g.\ \cite{Cas}, \cite{FV}, \cite{Haz}, \cite{Neu}, \cite{Serre}, and \cite{Serre2}.

\paragraph{{\bf Acknowledgments}} The author thanks his former fellow students at Harvard University, especially Jay Pottharst, who read the first draft and gave valuable comments. The revision of this paper was helped by the careful reading of Brian Conrad and suggestions by the referee. Ideas for simplification came through giving a course at the University of Cambridge in Fall 2007, and the author is grateful to all who attended the course. This work was partially supported by the EPSRC grant on Zeta Functions from the University of Nottingham, during the author's stay at Nottingham in the summer of 2005. The author was supported by the Society of Fellows at Harvard University and the Clay Mathematics Institute during the revision period.

\end{document}